\documentclass[reqno,11pt]{amsart}
\usepackage{amsmath, amssymb, amsthm,amsfonts} 
\usepackage[english]{babel}
\usepackage{bbm}
\usepackage{graphicx}
\usepackage{url}
\usepackage{epstopdf}
\usepackage[a4paper,bindingoffset=0.5cm,left=2cm,right=2cm,top=2.5cm,bottom=2cm,footskip=.8cm]{geometry}

\usepackage{tikz}
\usetikzlibrary{arrows, automata,positioning,calc,shapes,decorations.pathreplacing,decorations.markings,shapes.misc,petri,topaths}
\usepackage{tkz-berge}
  \usepackage{pgfplots}
  \pgfplotsset{compat=newest}
  \usetikzlibrary{plotmarks}
  \usepackage{grffile}
\newlength\figureheight
  \newlength\figurewidth
\pgfplotsset{%
    tick label style={font=\scriptsize},
    label style={font=\footnotesize},
    legend style={font=\footnotesize},
         every axis plot/.append style={very thick}
}

\usepackage{rotating}
\usepackage{amsbsy,enumerate}
\usepackage{graphicx}
\usepackage{ccaption}
\usepackage{comment}
\usepackage{mathrsfs}

\newcommand{\bs}{\boldsymbol}

\newcommand{\vb}{\vspace{3.2mm}}

\allowdisplaybreaks

\newtheorem{remark}{Remark}

\newtheorem{proposition}{Proposition}

\makeatletter
\renewcommand{\fnum@figure}[1]{\textbf{\figurename~\thefigure}. }
\renewcommand{\fnum@table}[1]{\textbf{\tablename~\thetable}. }
\makeatother

\begin{document}

\title[Queues on a dynamically evolving graph]{Queues on a dynamically evolving graph}

\author{Michel Mandjes, Nicos Starreveld, Ren\'e Bekker}

\begin{abstract}
This paper considers a population process on a dynamically evolving graph, which can be alternatively interpreted as a queueing network.
The queues are of infinite-server type, entailing that at each node all customers present are served in parallel. The links that connect the queues have the special feature that they are {\it unreliable},
in the sense that their status alternates between `up' and `down'. 
If a link between two nodes is down,
with a fixed probability each of the clients attempting to use that link is lost; otherwise the client remains at the origin node and reattempts using the link (and jumps
to the destination node when it finds the link restored).  For these networks we present the following results:
(a)~a system of coupled partial differential equations that describes the joint probability generating function 
corresponding to the queues' time-dependent behavior (and a system of ordinary differential equations for its stationary counterpart), (b)~an
algorithm to evaluate the (time-dependent and stationary) moments, and procedures to compute user-perceived performance measures which facilitate the quantification of the impact of the links' outages, (c)~a diffusion limit for the joint queue length process. We include explicit results for a series relevant special cases, such as tandem networks and symmetric fully connected networks.

\vb


\vb

\noindent
{\sc Affiliations.} M.~Mandjes and N.\ Starreveld are with Korteweg-de Vries Institute, University of Amsterdam.
R.\ Bekker is with Department of Mathematics, Vrije Universiteit, Amsterdam. 
The research for this paper is partly funded by the NWO Gravitation Programme N{\sc etworks}, Grant Number 024.002.003 (Mandjes, Starreveld), and an NWO Top Grant, Grant Number 613.001.352 (Mandjes).

\end{abstract}

\maketitle

\section{Introduction}
When considering a population process on a graph, the underlying network is typically assumed to be {\it static}: the network structure (i.e., the set of links that connect the nodes) is assumed to be constant over time. In 
many real-life situations, however, links may be temporarily inactive, entailing that the underlying structure should instead be considered as {\it dynamic}.  At a conceptual level, such a system can be seen as a network of queues, where the links' availability fluctuates in time. The main objective of this paper is to study the performance of such a queueing network, with links alternating between being `up' and `down'. Leading examples in which our model can be used include communication networks, road traffic networks, various physics-motivated networks, and chemical reaction networks. 

At a somewhat more detailed level, the network can be described as follows. The network is a graph with nodes and links, along which clients travel. At any node, external arrivals occur according to a Poisson process with a node-specific rate. Service times at the nodes are exponentially distributed (with a node-specific parameter); when a customer has been served at a node, he selects a next node through some routing mechanism (where it is also an option to leave the network). Suppose the client resides at node $i$ and he wants to be routed to node $j$; assume that each link's up- and down-times are exponentially distributed. Then, depending on the situation at hand, the following two options arise. If the link from $i$ to $j$ is up, then he jumps from node $i$ to node $j$. If, on the contrary,  the link from $i$ to $j$ is down, then either  the client is lost (which happens with a node-specific probability), or he  waits an exponentially distributed amount of time at node $i$ and tries again. 

The queueing mechanism studied in this paper is infinite-server, making our analysis particularly useful for situations in which there is no (or hardly any) interference between the clients at each individual queue, in the sense that they can be served essentially in parallel. It is noted that in this paper we use queueing-theoretic terminology, but infinite-server queues are frequently used in other domains as well. As a model in which particles move on a dynamically evolving graph, it can be seen as an object relevant to statistical physics (cf., for instance, the model considered in 
\cite{DH}), but there are applications in chemical reaction networks \cite{KUR}, (cell) biology \cite{BRUG}, and population dynamics \cite{LOU} as well. In operations research, the infinite-server model we have defined can be used to study e.g.\ the  numbers of clients simultaneously using (somewhat larger) segments in a road traffic network, or numbers of clients simultaneously visiting connected websites. 

For this class of model, we are interested in various performance measures. The most important one is the joint distribution of the (time-dependent and stationary)  queue lengths at all nodes, together with the number of 
{\it lost} clients (i.e., clients who leave the network because of  link failures). 
\begin{itemize}
\item[$\circ$]
In the first part of the paper we derive a series of exact results. 
(i)~Our first class of results is in terms of a system of coupled partial differential equations for the probability generating function pertaining  to the joint queue length distribution. (ii)~In the second place, this system of differential equations can be used to recursively determine all (time-dependent and stationary) moments. (iii)~Thirdly, we assess the impact of the network's down-times on the service quality that is perceived by its users. 
\item[$\circ$]
Then we consider scaling limits: by scaling the external arrival rates and the up- and down-times, we present a diffusion limit. This result entails that the joint queue-length process weakly converges to a mean-reverting Gaussian process (viz.\ a multivariate Ornstein-Uhlenbeck process). An important feature of the scaling chosen is that the speed at which the external arrival rates are scaled may differ from the speed at which the up- and down-times are scaled. This  creates the flexibility to cover networks in which the alternation between up- and down times is relatively slow (think of road networks) or relatively fast (think of the channel conditions in a wireless network); also time-scale separation ideas (as often relied on in chemical reaction networks) can  thus be modeled.
\end{itemize}

The model in this paper can be seen as an instance of a stochastic process (viz.\ a queueing process) on a dynamically evolving graph. The literature on such models is still at its infancy. Where static random graphs form a classical topic in probability theory, dating back to the pioneering work of Erd\H{o}s and R\'enyi \cite{ER} and Gilbert \cite{EG}, only recently the behavior of randomly evolving graphs has received substantial attention; see e.g.\ \cite{HO1,HO2,MSBS,ZMN} for a few examples. Examples of papers on random processes on (dynamic) random graphs are \cite{AV, BM, BR1}.
The systematic study of queueing processes on such a randomly evolving graph has hardly been looked at, a notable exception being the recent study \cite{FMP}.  The model considered in \cite{FMP} complements the one studied in the present work. Most notably, the framework of \cite{FMP} in particular facilitates modelling the effect of nodes going down every now and then, where the present paper has a focus on links going down. The immediate consequence of this difference in modelling, is that in the framework of \cite{FMP} diffusion limits do not apply, due to  the instantaneous downward jumps of the network population vector at epochs that a node fails. 

There is also a relation with the classical work
%by Massey and Whitt 
\cite{MW}, where the Poisson-arrival-location model (PALM) is introduced. In this model customers arrive according to an inhomogeneous Poisson process and move independently through the network according to some random location process (with a fixed routing matrix). A consequence of the way the model is constructed is that, for instance, the number of customers at each node follows a Poisson distribution. The major difference with our model, is that in our setup the topology of the network is determined by a modulating process (meaning that the routing matrix is random); consequently the positions of different clients (during their path through the network) are in our model no longer independent, thereby also destroying the `Poisson properties'. Observe that  it is this dependence structure that considerably complicates the analysis. It also explains why we pursue scaling limits for obtaining insight in the network population distribution (which is obviously not needed in the setup of \cite{MW}, as there closed-form expressions are available).

Our analysis will be based on casting our model as a network of infinite-server queues under Markov modulated arrival and service rates. Explicit results on (single-node) Markov modulated infinite-server queues (primarily in terms of differential equations for the probability generating function, and the corresponding moments)
 can be found in e.g.\ \cite{DAU,GUR,FRA,KEI,OCP}. Diffusion limits for such single-node systems have been derived in e.g.\ \cite{AND,BKM}; we also refer to \cite{JMTW} for a recent contribution with such diffusion results for a broad class of networks of Markov modulated infinite-server queues. For general background on queueing networks, we refer to \cite{KEL,KR,SER}.

\vb

This paper is organized as follows. In Section 2 we describe our model. Section 3 presents our analysis, in terms of results exact results for the probability generating function and moments; we restrict ourselves to the case that clients who wish to jump but the corresponding link is down, are lost with probability 1. 
In Section 4 concerns the weak convergence to a Gaussian process, for the same model. In Section 5 we consider a number of extensions, including the one in which blocked customers are not necessarily lost but retry. In Section 6 we discuss a number of special cases for which the calculations can be done explicitly. Concluding remarks are found in Section 7. 

\section{Model description}
In this section we first provide a detailed model description, and then introduce   quantities of our interest. 

\vb

The network that we consider consists of $n$ nodes that are connected through $\bar n:={ {n}\choose{2}}$ links. 
Let $\lambda_i$ be the rate of the Poissonian arrival process at node $i$. The time spent at node $i$ is exponentially distributed with parameter $\mu_i$ (where we discuss in Section 5 how our setup extends to the case of phase-type service times). After having been served at node $i$, the probability that the served customer wishes to jump to node $j$ (where $j\not=i$) is $p_{ij}$, where $p_{i0}$ is the probability of leaving the network. We obviously assume that $\sum_{j\not=i}p_{ij}=1$, and we write $\mu_{ij}=\mu_i p_{ij}$.
It is noted that  this setup does {\it not} necessarily mean that we assume that the network be a complete graph; if a node pair $(i,j)$ is not connected, we are to set the corresponding $\mu_{ij}$ equal to $0$.
Observe that the dynamics as described above entail that the number of clients evolves as  an infinite-server queue: the clients are served in parallel, and hence do not interact. 
We assume that the routing mechanism gives rise to an irreducible structure, entailing that the $\mu_{ij}$ are such that for a client residing at a specific node with positive probability it visits any other node before leaving the network. In addition, for at least one node $i$ it holds that $\mu_{i0}$ is {\it strictly} positive, thus guaranteeing that the network is stable. The arrival processes and service/routing processes are assumed independent.

\vb

We now describe how the links alternate between being `up' and `down'. To this end, we let the underlying graph dynamics be determined by a $K$-dimensional background process $({\bs X}(t))_{t\geqslant 0}$,  assumed to be independent of the arrival processes and the serving/routing processes,  that is defined as follows. The $\bar n$ links are  partitioned into  $K$ mutually disjoint sets, which are denoted by $A_1, \ldots , A_K$, which we refer to as {\it blocks}. All links that lie in a specific block, say $A_k$,
alternate between `up' and `down' simultaneously; $X_k(t)=1$ means that at time $t$ the links in block $k$ are `up', and $0$ otherwise. We define
\[Q^{(k)}=\left(\begin{array}{rr}-q^{(k)}_0&q^{(k)}_0\\
q^{(k)}_1&-q^{(k)}_1\end{array}\right);\]
the down-time (up-time, respectively) of block $k$ is exponentially distributed with parameter $q^{(k)}_0$
($q^{(k)}_1$). The `graph process' is given through
\[({\bs X}(t) )_{t\geqslant 0}= (X_1(t),\ldots, X_K(t))_{t\geqslant 0},\]
which attains values in $\{0,1\}^K.$
The two extreme scenarios are on one hand the case that we have just one block consisting of all $\bar n$ links, or on the {other}  hand the case that we have $\bar n$ independently evolving blocks that consist of one link each.  The transition rate matrix  of ${\bs X}(\cdot)$ is of dimension $\bar K\times \bar K$ with $\bar K:=2^K$, and given by
\[{\bs Q}:=\bigoplus_{k=1}^KQ^{(k)}=\sum_{k=1}^K I_{2^{k-1}}\otimes Q^{(k)}\otimes I_{2^{K-k}},\]
where $I_n$ denotes the $n$-dimensional identity matrix and $B_1 \otimes B_2$ denotes the Kronecker product of the two matrices $B_1$ and $B_2$. We let $q_{k\ell}$ be the $(k,\ell)$-th entry of ${\bs Q}$.

We now explain what happens to a client who wants to jump from $i$ to $j$ when the link is not present. As long as the link is down, at any attempt the client is lost with probability $f_{ij}\in[0,1]$, and he remains at the node with probability $1-f_{ij}$. While being at the node the mechanism that we defined above is in place: after an exponentially distributed amount of time with mean $\mu_{ij}^{-1}$ (with $j=0,\ldots,n$) he wishes to jump to node $j$. To keep the notation compact, in Sections 3 and 4 we assume that $f_{ij}=1$ for all $i,j=1,\ldots,n$ (i.e., all clients are lost who wish  to jump from $i$ to $j$ when the link between $i$ and $j$ is absent); in Section 5 we point out how to adapt the results to include situations with $f_{ij}\in[0,1).$

\vb

In this paper a key role is played by the $n$-dimensional queue length process
\[({\bs M}(t))_{t\geqslant 0} = (M_1(t),\ldots, M_n(t))_{t\geqslant 0},\]
where $M_i(t)\in{\mathbb N}_0$ represents the number of clients in the queue at node $i$ at time $t$. Our objective is to characterize the distribution of ${\bs M}(t)$; as we will see below, this is possible, albeit in implicit terms, viz.\ in terms of a partial differential equation for the corresponding joint probability generating function.  Observe that by itself $({\bs M}(t))_{t\geqslant 0}$ is not a Markov process, but the joint process $({\bs M}(t),{\bs X}(t))_{t\geqslant 0}$ is. 

As we want  to keep track of ${\bs M}(t)$ as well as the number of lost clients,
we work with the probability generating function
\[{\varphi}_k(w,{\bs z},t) = {\mathbb E}\, \big[w^{L(t)}z_1^{M_1(t)}\cdots z_n^{M_n(t)} 1_{\{{\bs X}(t) = k\}}\big],\]
with $L(t)$ defined as the number of lost clients  {due to a link being `down'} during the interval $[0,t]$
and  $k$ be an element in $\{0,1\}^K$. 

\begin{remark}{\em
In the model described all links are bidirectional: if the link between $i$ and $j$ {is `down'}, then clients can jump neither form $i$ to $j$ nor from $j$ to $i$. The unidirectional variant of our model works in the precise same way; then there are $n(n-1)$ (instead of $\frac{1}{2}n(n-1)$) possible links.}
\hfill$\Diamond$
\end{remark}

\section{Prelimit results}
In this section we first set up a system of coupled partial differential equations for ${\bs \varphi}(w,{\bs z},t)$ (i.e., the $2^K$-dimensional vector with elements ${ \varphi}_{k}(w,{\bs z},t)$).
We then point out how these can be used to determine moments. The next subsection presents ways to quantify the effect of the graph dynamics on the performance as perceived by the network's users.

\subsection{Partial differential equations}
The main idea is to express ${\bs \varphi}(w,{\bs z},t+\Delta t)$, for $\Delta t$ small, in terms of ${\bs \varphi}(w,{\bs z},t)$. We follow the precise same procedure as in e.g.\ \cite{KOEN}: we first set up the Kolmogorov equations for the state $(L(t),{\bs M}(t)) = (m_0,\ldots,m_n)$ at time $t\geqslant 0$, then multiply with $w^{m_0}z_1^{m_1}\cdots z_n^{m_n}$, and sum over $m_0,\ldots,m_n$. We recognize probability generating functions and their derivatives. 
More specifically, with ${\mathbb I}(i,j,k)$ being 1 if the link $(i,j)$ is `up' when ${\bs X}(\cdot)$ is in state $k$ and $0$ {otherwise}, we thus obtain,
\begin{align*}
{ \varphi}_{k}(w,{\bs z},t+\Delta t)=\,&{\varphi}_{k}(w,{\bs z},t)+
\sum_{i=1}^n{ \varphi}_{k}(w,{\bs z},t)(z_i-1)\cdot\lambda_i \,\Delta t\,+\\
&\sum_{i=1}^n\sum_{ {j=1,}j\not=i}^n\frac{\partial {  \varphi}_{k}(w,{\bs z},t)}{\partial z_i}\left( {z_j}-{z_i}\right)
\cdot {\mathbb I}(i,j,k)\cdot\mu_{ij} \,\Delta t \,+
\\&\sum_{i=1}^n\sum_{{j=1,}j\not=i}^n\frac{\partial {  \varphi}_{k}(w,{\bs z},t)}{\partial z_i}\left( w-{z_i}\right)
\cdot \big(1-{\mathbb I}(i,j,k)\big)\cdot\mu_{ij} \,\Delta t \,+
\\&\sum_{i=1}^n\frac{\partial {  \varphi}_{k}(w,{\bs z},t)}{\partial z_i}\left(1-{z_i}\right)\cdot\mu_{i0} \,\Delta t\,+\\
&\sum_{\ell\not=k}{ \varphi}_{\ell}(w,{\bs z},t)\cdot q_{\ell k}\,\Delta t-
\sum_{\ell\not=k}{  \varphi}_{k}(w,{\bs z},t)\cdot q_{k\ell}\,\Delta t+o(\Delta t).
\end{align*}
The next step is to subtract ${\varphi}_{k}(w,{\bs z},t)$ from both sides, divide by $\Delta t$, and send $\Delta t\downarrow 0$. In matrix-vector form, the resulting system of coupled partial differential equations reads, with ${\mathbb I}_{\bar K}(i,j):={\rm diag}\{{\mathbb I}(i,j,1),\ldots,{\mathbb I}(i,j,\bar K)\}$
and ${\mathbb J}_{\bar K}(i,j):=I_{\bar K}-{\mathbb I}_{\bar K}(i,j)$, as follows. The function ${\bs  \varphi}({\bs z})$ denotes the probability generating function of the stationary counterpart ${\bs M}$ of $({\bs M}(t))_{t\geqslant 0}$.
\begin{proposition} \label{pr} The joint probability generating function ${\bs  \varphi}(w,{\bs z},t)$ satisfies
\begin{align*}\frac{\partial {\bs  \varphi}(w,{\bs z},t)}{\partial t} = &\,
\sum_{i=1}^n{\bs \varphi} (w,{\bs z},t)\,\lambda_i(z_i-1) \,+\sum_{i=1}^n\sum_{{j=1,}j\not=i}^n\frac{\partial {\bs  \varphi}(w,{\bs z},t)}{\partial z_i}
{\mathbb I}_{\bar K}(i,j)\,
\mu_{ij}\left( {z_j}-{z_i}\right)\,+
 \\
&\sum_{i=1}^n\sum_{{j=1,}j\not=i}^n\frac{\partial {\bs  \varphi}(w,{\bs z},t)}{\partial z_i}{\mathbb J}_{\bar K}(i,j)\,
\mu_{ij}\left( w-{z_i}\right)\,+\\
&\sum_{i=1}^n \frac{\partial {\bs  \varphi}(w,{\bs z},t)}{\partial z_i}\,\mu_{i0} (1-z_i)
+
 {\bs  \varphi}(w,{\bs z},t)\, {\bs Q}
.\end{align*}
The probability generating function of the stationary counterpart 
${\bs  \varphi}({\bs z})$ satisfies
\begin{align*} {\bs  0} = &\,
\sum_{i=1}^n{\bs \varphi} ({\bs z})\,\lambda_i(z_i-1) \,+\sum_{i=1}^n\sum_{j=1,j\not=i}^n\frac{\partial {\bs  \varphi}({\bs z})}{\partial z_i}
{\mathbb I}_{\bar K}(i,j)\,
\mu_{ij}\left( {z_j}-{z_i}\right)\,+
 \\
&\sum_{i=1}^n\sum_{j=1,j\not=i}^n\frac{\partial {\bs  \varphi}({\bs z})}{\partial z_i}{\mathbb J}_{\bar K}(i,j)\,
\mu_{ij}\left( 1-{z_i}\right)\,+\\
&\sum_{i=1}^n \frac{\partial {\bs  \varphi}({\bs z})}{\partial z_i}\,\mu_{i0} (1-z_i)
+
 {\bs  \varphi}({\bs z})\, {\bs Q}
.\end{align*}
\end{proposition}

\subsection{First moment} 
\label{M1}
In this section we exploit Proposition \ref{pr} to determine the first moments; we first point out how this procedure works for the stationary queue length ${\bs M}$, but later indicate how the corresponding transient moments can be found. We let ${\bs X}$ denote the stationary version of the background process.

Define, for $i=1,\ldots,n$,
\[{\bs v}_i:= ({\mathbb E} M_i 1_{\{{\bs X}=1\}},\ldots,{\mathbb E} M_i 1_{\{{\bs X}=\bar K\}})=
\lim_{{\bs z}\uparrow {\bs 1}}
\frac{\partial {\bs  \varphi}({\bs z})}{\partial z_i}
.\]
Let ${\bs \pi}$ be the invariant probability measure of ${\bs Q}$, i.e., the $\bar K$-dimensional row-vector such that ${\bs \pi}{\bs Q}={\bs 0}$ and whose entries sum to 1. 
By differentiating the differential equation featuring in Proposition \ref{pr} with respect to $z_i$ and letting   ${\bs z}\uparrow {\bs 1}$, we obtain, for $i=1,\ldots,n$,
\[{\bs 0} = {\bs \pi}\lambda_i -\sum_{j=1,j\not=i}^n{\bs v}_i \,\mu_{ij}+
\sum_{j=1,j\not=i}^n{\bs v}_j \,\mu_{ji}{\mathbb I}_{\bar K}(j,i)-{\bs v}_i\,\mu_{i0}+{\bs v}_i\,{\bs Q}
.\]

We now explain how to set up a computational procedure with which the ${\bs v}_i$ can be found. The $n$ sets of $\bar K$-dimensional systems of linear  equations can be cast into a single set of $n\bar K$ linear equations (in equally many unknowns). Let ${\bs v}\equiv ({\bs v}_1,\ldots,{\bs v}_n)$.
Also, let ${\bs \lambda}$ the row-vector $(\lambda_1,\ldots,\lambda_n)$, and
\[\nu_i:=\sum_{j=1,j\not=i}^n \mu_{ij}.\] 
In addition, we define the matrices
\[{\mathscr M}_+:=\left(\begin{array}{cccc}
\nu_1 I_{\bar K}&&&\\
&\nu_2 I_{\bar K}&&\\
&&\ddots&\\
&&&\nu_n I_{\bar K}\end{array}\right)\hspace{-0.8mm},\:\:
{\mathscr M}_-:=\left(\begin{array}{cccc}
0 &\mu_{12}{\mathbb I}_{\bar K}(1,2)&\hdots&\mu_{1n}{\mathbb I}_{\bar K}(1,n)\\
\mu_{21}{\mathbb I}_{\bar K}(1,2)&0&&\mu_{2n}{\mathbb I}_{\bar K}(2,n)\\
\vdots&&\ddots&\\
\mu_{n1}{\mathbb I}_{\bar K}(n,1)&\mu_{n2}{\mathbb I}_{\bar K}(n,2)&&0\end{array}\right)\hspace{-0.8mm},\]
and
\[{\mathscr M}_0:=\left(\begin{array}{cccc}
\mu_{10} I_{\bar K}&&&\\
&\mu_{20} I_{\bar K}&&\\
&&\ddots&\\
&&&\mu_{n0} I_{\bar K}\end{array}\right)\hspace{-0.8mm},\:\:
{\mathscr Q}:=\left(\begin{array}{cccc}
{\bs Q}&&&\\
&{\bs Q}&&\\
&&\ddots&\\
&&&{\bs Q}\end{array}\right)\hspace{-0.8mm}.\]
We thus arrive at the linear system
\[{\bs\lambda} \otimes {\bs \pi} ={\bs v}( {\mathscr M}_+-{\mathscr M}_-+{\mathscr M}_0-{\mathscr Q}),\]
so that ${\bs v}= ({\bs\lambda} \otimes {\bs \pi}) {\mathscr N}^{-1},$ with ${\mathscr N}:=
{\mathscr M}_+-{\mathscr M}_-+{\mathscr M}_0-{\mathscr Q}.$

The transient first moment follows immediately by solving the corresponding system of linear differential equations. Let $({\bs M}(0), X(0))=({\bs m},k_0)$, and let ${\bs e}_k$ the $k$-th unit vector. 
Then, in self-evident notation, and with ${\bs \pi}(t)={\bs e}_{k_0}^{\rm T} \,{\rm e}^{{\bs Q}t}$,
\[{\bs v}'(t) = {\bs\lambda} \otimes {\bs \pi}(t) - {\bs v}(t) {\mathscr N},\]
which is solved by
\begin{equation}\label{MEAN}{\bs v}(t) ={\bs m} \,{\rm e}^{-{\mathscr N}t} +\int_0^t ({\bs\lambda} \otimes {\bs \pi}(s))
\, {\rm e}^{-{\mathscr N}(t-s)}{\rm d}s.\end{equation}
This approach can be extended in a straightforward fashion to also include the mean of the number of clients lost.

\begin{remark}{\em The relation (\ref{MEAN}) can be regarded as the expected-value version of the distributional identity
\[{\bs M}(t) \stackrel{\rm d}{=}{\bs m} \,{\rm e}^{-{\mathscr N}t} +\int_0^t ({\bs\lambda} \otimes {\bs X}(s))
\, {\rm e}^{-{\mathscr N}(t-s)}{\rm d}s,\]
cf.\ the relation used for our model's single-queue counterpart in 
\cite{BKMT}. The above distributional equality has an insightful interpretation: the first term on the right-hand side 
corresponds to the contribution to ${\bs M}(t)$ of clients that were already present at time $0$, whereas the second term 
represents the contribution of arrivals in $[0,t]$ which are then appropriately `thinned'.

As pointed out in \cite{BKMT}, this representation in principle also provides a method 
to evaluate the covariance matrix of ${\bs M}(t)$, using the law of total variance; in Section \ref{HM} we will rely on an alternative approach, though. 
\hfill$\Diamond$
}\end{remark}

\begin{remark} \label{rem1}{\em In the fully symmetric situation, the formulas simplify considerably. Let $\lambda$ be the arrival rate at each of the $n$ nodes. The service rate is $\sigma:=\nu+\mu_0$, where the client leaves the network with probability $\mu_0/\sigma$ and wants to move to another node (which is then picked uniformly at random) with probability $\nu/\sigma.$ Let all blocks alternate independently between being `up' and `down'. The up- and down rates are denoted by $q_0$ and $q_1$, respectively (i.e., the down-time is exponentially distributed with mean $q_0^{-1}$, and the up-time exponentially distributed with mean $q_1^{-1}$).
%, such that $q_0/q_1$ is constant across all blocks. 
Assume that the queues start empty at time $0$, while  all links are in stationary state (i.e., each of them is  `up' with
probability $\pi:= q_0/(q_0+q_1)$). 
Then, for each of the {nodes} the mean number of clients present satisfies
\[v'(t) = \lambda + (n-1)\,v(t) \frac{\nu\pi}{n-1} - v(t)\sigma= \lambda-(\nu(1-\pi)+\mu_0)v(t).\] 
We thus find that
\[v(t) =\frac{\lambda}{\nu(1-\pi)+\mu_0} \left(1-{\rm e}^{-(\nu(1-\pi)+\mu_0)t}\right),\]
which converges to $v:= \
{\lambda}/({\nu(1-\pi)+\mu_0})$ 
as $t\to\infty$. \hfill$\Diamond$
}\end{remark}

\subsection{Higher moments}\label{HM}
We now point out how (mixed) higher moments can be evaluated. We work here, for obvious reasons, with the {\it factorial} moments, from which the regular moments can be recovered in an evident manner.
We recall the standard notation $(m)_r:=m!/(m-r)!$ for $m,r\in{\mathbb N}$ with $m\geqslant r$; this notation is the so-called {\it Pochhammer symbol}.  The objective here is to compute, for ${\bs r}\equiv (r_1,\ldots,r_n)$ with $r_i\in{\mathbb N}_0$,
\[\psi_k({\bs r},t):={\mathbb E}\left(\left(\prod_{i=1}^{n}(M_i(t))_{r_i}\right) 1_{\{{\bs X}(t)=k\}}\right) = \lim_{w\uparrow 1,{\bs z}\uparrow {\bs 1}}
\frac{\partial^{r_1+\cdots+r_n}\varphi_k(w,{\bs z},t)}{\partial z_1^{r_1}\cdots \partial z_n^{r_n}}.\] 
We will show that the $\psi_k({\bs r},t)$ can be recursively evaluated. To this end, we first introduce the following differential operator: for $f: {\mathbb R}\times {\mathbb R}^n\times {\mathbb R}\mapsto {\mathbb R}$,
\[{\mathbb D}_{\bs r}[f(w,{\bs z},t)]:= \frac{\partial^{r_1+\cdots+r_n}f(w,{\bs z},t)}{\partial z_1^{r_1}\cdots \partial z_n^{r_n}}.\] 
Now the idea is to impose the operator ${\mathbb D}_{\bs r}[\cdot]$ on both sides of the partial differential equation given in Proposition \ref{pr}. In addition considering the limit $w\uparrow 1,{\bs z}\uparrow {\bs 1}$, we thus obtain, with ${\bs e}_i$ the $i$-th $n$-dimensional unit vector, and $\mu_{ijk}^+:= {\mathbb I}(i,j,k)\,\mu_{ij}$ and $\mu_{ijk}^-:= (1-{\mathbb I}(i,j,k))\,\mu_{ij}$,
%\footnote{\tt\tiny MM: this computation is not entirely trivial, with lots of indices floating around, so I may have made mistakes. Please check!}
\begin{align*}
\lim_{w\uparrow 1,{\bs z}\uparrow {\bs 1}}{\mathbb D}_{\bs r}[ \frac{\partial {\varphi}_k(w,{\bs z},t)}{\partial t} ]&=\lim_{w\uparrow 1,{\bs z}\uparrow {\bs 1}}\Bigg(\sum_{i=1}^nr_i
{\mathbb D}_{{\bs r}-{\bs e}_i}[{\varphi}_k(w,{\bs z},t)]\lambda_i \,1_{\{r_i\not= 0\}}
\\
&\hspace{6mm}
+\sum_{i=1}^n\sum_{j=1,j\not=i}^n r_j {\mathbb D}_{{\bs r}-{\bs e}_j+{\bs e}_i} [{\varphi}_k(w,{\bs z},t)]
\mu_{ijk}^+ \,1_{\{r_j\not= 0\}}+ \\
&\hspace{6mm}-
\sum_{i=1}^n\sum_{j=1,j\not=i}^n r_i {\mathbb D}_{{\bs r}} [{\varphi}_k(w,{\bs z},t)]\mu_{ijk}^+-
\sum_{i=1}^n\sum_{j=1,j\not=i}^n r_i {\mathbb D}_{{\bs r}} [{\varphi}_k(w,{\bs z},t)]\mu_{ijk}^-\\
&\hspace{6mm}
-\sum_{i=1}^n r_i {\mathbb D}_{\bs r}[ {\varphi}_k(w,{\bs z},t)] \mu_{i0}+\sum_{\ell=1}^{\bar K}{\mathbb D}_{\bs r}[ {\varphi}_\ell(w,{\bs z},t)] q_{\ell k} \Bigg);
\end{align*}
these computations rely on the evident relation
\[\lim_{z\uparrow 1} \frac{\partial^r}{\partial z^r}f(z)(z-1) = \lim_{z\uparrow 1}rf^{(r-1)}(z)= rf^{(r-1)}(1).\]
Observe that in the above differential equation the indicator functions $1_{\{r_i\not= 0\}}$ (first term on the right-hand side) and $1_{\{r_j\not= 0\}}$ (second term on the right-hand side) can be left out: if the indicator function is $0$, the corresponding term equals $0$ anyway.

Consequently, the transient mixed reduced moments can be alternatively expressed in compact notation as follows. Here it is used  that $\mu_{ijk}^++\mu_{ijk}^-=\mu_{ij}$.
\begin{proposition} \label{P2}
For ${\bs r}\in {\mathbb N}_0^n$, $k\in\{1,\ldots,\bar K\}$ and $t\geqslant 0$,
\begin{align}\nonumber  \frac{\partial \psi_k({\bs r},t)}{\partial t}=\,&\sum_{i=1}^nr_i
\psi_k({\bs r}-{\bs e}_i,t)\lambda_i+ \sum_{i=1}^n\sum_{j=1,j\not=i}^n r_j  \psi_k({\bs r}-{\bs e}_j+{\bs e}_i,t)
\mu_{ijk}^+\,-\\
\,&\sum_{i=1}^n\sum_{j=1,j\not=i}^n r_i \psi_k({{\bs r}},t)\,
\mu_{ij}- \sum_{i=1}^n  r_i \psi_k({{\bs r}},t)\,
\mu_{i0}+ \sum_{\ell=1}^{\bar K} \psi_\ell({\bs r},t) q_{\ell k}.\label{dve2}
\end{align}
\end{proposition}
To obtain the stationary mixed reduced moments, one has to set the left-hand side equal to $0$.

The reduced moments can be determined recursively, by solving a  non-homogeneous system of linear differential equations. To verify this claim, define $\xi({\bs r})$ as the sum of the entries of ${\bs r}$, i.e., $r_1+\cdots +r_n.$ Let ${\mathscr S}_r$ be all vectors ${\bs r}$ such that $\xi({\bs r})=r.$ Observe that 
\[\xi({\bs r}) = \xi({\bs r}-{\bs e}_j+{\bs e}_i) = \xi({\bs r}-{\bs e}_i)+1.\]
The idea now is to use {(\ref{dve2})} to evaluate $( \psi_1({\bs r},t),\ldots, \psi_{\bar K}({\bs r},t))$ recursively: for ${\bs r}$ such that $\xi({\bs r})=r$ this vector is computed using the corresponding expressions for ${\bs r}$ such that $\xi({\bs r})=r-1$. In more detail, this procedure works as follows.
\begin{itemize}
\item[$\circ$]
For $r=0$, we find the $\psi_k({\bs 0},t)={\mathbb P}(X(t)=k)= ({\rm e}^{{\bs Q}t})_{k_0,k}.$
\item[$\circ$]
For $r=1$, we find the $\psi_k({\bs e}_i,t)$ (for $i=1,\ldots,n$) by appealing to Eqn.\ {(\ref{dve2})}, 
using the $\psi_k({\bs 0},t)$ that we found for $r=0$. This amounts to solving $\bar K n$ coupled linear differential equations (where it can be checked that this system is equivalent to the one set up in Section \ref{M1}).
\item[$\circ$]
We now consider $r=2$. It is readily checked that $\#\{{\mathscr S}_2\} = n+\frac{1}{2}n(n-1)=\frac{1}{2}n(n+1)$, and that there are equally many equations of the type (\ref{dve2}). As a consequence, using the $\psi_k({\bs e}_i,t)$ that we found for $r=1$, Eqn.\ (\ref{dve2}) reveals that we have to solve a non-homogeneous system of $\bar K\cdot\frac{1}{2}n(n+1)$ linear differential equations. 
\item[$\circ$]
One can proceed in a similar way with $r\geqslant 3$; e.g. for $\xi({\bs r})=3$ we have that 
\[\#\{{\mathscr S}_3 \}= n + 2\cdot\textstyle\frac{1}{2} n(n-1) +\textstyle\frac{1}{6}n(n-1)(n-2)= \textstyle\frac{1}{6}n(n+1)(n+2).\]
The cases with $\xi({\bs r})\geqslant 3$ are solved very similarly to the case $\xi({\bs r})=2$. Using (\ref{dve2}) it can be inductively verified that in the $r$-th step we have a system of $K_r:=\bar K\cdot \#\{{\mathscr S}_r \} $
non-homogeneous  linear differential equations, where
\[\#\{{\mathscr S}_r \} ={{n+r-1}\choose{r}}= \frac{1}{r!}\cdot  n(n+1)\cdots(n+r-1).\] 
\end{itemize}
For the stationary reduced means a similar recursive procedure can be set up; in the $r$-th step a $K_r$-dimensional system of linear equations needs to be solved. 

The above procedure can be extended in an evident way to include also the number of customers lost, focusing on the object
\[\bar\psi_k({\bs r},t):={\mathbb E}\left(\left((L(t))_{r_0}\prod_{i=1}^{n}(M_i(t))_{r_i}\right) 1_{\{{\bs X}(t)=k\}}\right) = \lim_{w\uparrow 1,{\bs z}\uparrow {\bs 1}}
\frac{\partial^{r_0+\cdots+r_n}\varphi_k(w,{\bs z},t)}{\partial w^{r_0}\partial z_1^{r_1}\cdots \partial z_n^{r_n}},\] 
with ${\bs r}\equiv (r_0,\ldots,r_n).$

\subsection{User-perceived performance}
In this subsection we study the impact of the links' outages on the performance as perceived by the users. 
As it turns out, this can be done relying on classical arguments. We start by analyzing the fraction of clients that are lost (i.e., clients who leave the network because of a link being down, and not because of a service completion), denoted by $\omega$. To this end, let $\omega_{ik}$ be the probability that a client who enters the network at node $i$ while the background process is in state $k$, is lost. 
Observe that ${\mu_{ijk}^+}/{\sigma_{ik}}$ is the probability of a client jumping from node $i$ to node $j$ when the background process is in $k$; ${q_{k\ell}} /{\sigma_{ik}}$ and ${\mu_{ijk}^-}/{\sigma_{ik}}$
can be interpreted analogously. 
Then, with $q_k:=-q_{kk}$ and $\sigma_{ik}:=\nu_i+\mu_{i0}+q_k$, by conditioning on the first jump,
\[\omega_{ik} = \sum_{j\not=i} \left(\frac{\mu_{ijk}^+}{\sigma_{ik}}\right) \omega_{jk}
+\sum_{\ell\not=k} \left(\frac{q_{k\ell}} {\sigma_{ik}}\right) \omega_{i\ell}+
\sum_{j\not=i} \left(\frac{\mu_{ijk}^-}{\sigma_{ik}}\right)
,\]
or, in a more convenient form,
\[-\sum_{j\not=i} {\mu_{ijk}^-} = -\sigma_{ik}\omega_{ik} +
\sum_{j\not=i} {\mu_{ijk}^+} \omega_{jk}
+\sum_{\ell\not=k}  q_{k\ell}  \omega_{i\ell}
.\]
These equations constitute an $n\bar K$-dimensional diagonally dominant system of linear equations (actually even {\it strictly} diagonally dominant as there is at least one $i$ such that $\mu_{i0}>0$), which is known to yield a unique solution. Let, as before, $\pi_\ell$ denote the stationary probability that the background process ${\bs X}(\cdot)$ is in state $\ell$ (i.e., ${\bs\pi}$ solves ${\bs\pi}{\bs Q} = {\bs 0}$). 
Then, with $\bar{\lambda} := \sum_{i=1}^n \lambda_i$, the loss probability equals
\[\omega = \sum_{k=1}^{\bar K} \pi_k\left(\frac{1}{\bar\lambda} {\sum_{i=1}^n \lambda_i \omega_{ik}}\right).\]
As an aside we mention that the loss probability $\omega$ can alternative be evaluated, using the methodology of the previous subsections, as
\[\lim_{t\to\infty}\frac{{\mathbb E}\,L(t)}{\bar\lambda t}.\]
Along the same lines, we can determine the mean time  (to be denoted by $\tau$) the job remains in the network, jointly with the client being eventually lost (i.e., leaving the network because of a link failure).
%; analogously, we find the corresponding quantity jointly with the client leaving the network without being lost. 
Define  $\tau_{ik}$ to be this quantity for a client who enters the network at node $i$ while the background process is in state $k$.
Then, again conditioning on the first jump,
\[\tau_{ik} =
 \sum_{j\not=i} \left(\frac{\mu_{ijk}^+}{\sigma_{ik}}\right)\left(\frac{1}{\sigma_{ik}}+ \tau_{jk}\right)
+\sum_{\ell\not=k} \left(\frac{q_{k\ell}} {\sigma_{ik}}\right)\left(\frac{1}{\sigma_{ik}}+ \tau_{i\ell}\right)+
\sum_{j\not=i} \left(\frac{\mu_{ijk}^-}{\sigma_{ik}}\right)\left(\frac{1}{\sigma_{ik}}\right).\]
Again this system of linear equations is diagonally dominant. 
As above,
\[\tau = \sum_{k=1}^{\bar K} \pi_k\left(\frac{1}{\bar\lambda} {\sum_{i=1}^n \lambda_i \tau_{ik}}\right).\]

\section{Scaling limit: functional central limit theorem} 
In this section we study the system under a particular scaling, under which there is convergence to a Gaussian process (viz.\ a multivariate Ornstein-Uhlenbeck process). 
We  consider the following scaling: ${\bs\lambda}\mapsto N{\bs\lambda}$ and $q_i^{(k)}\mapsto N^\alpha q_i^{(k)}$ for some $\alpha>0$ (for $k=1,\ldots,K$ and $i=0,1$). The model is thus parametrized by $N$; to stress the dependence on $N$, we throughout write ${\bf M}^{(N)}(t)$ and ${\bs X}^{(N)}(t)$.  We focus on a functional central limit theorem for a {centered} and appropriately scaled version of ${\bf M}^{(N)}(t)$. For the moment we concentrate on the (more involved) case $\alpha=1$; in Remark~\ref{ralpha} we point out how this provides us with the limit result  for $\alpha\in(0,\infty)\setminus \{1\}$ as well.

A full proof is (far) beyond the scope of this paper. For a rigorous derivation of the functional central limit theorem, based on martingale arguments in combination with the continuous mapping theorem, for a class of network models that is substantially broader than the one studied in this paper, we refer to \cite{JMTW}. Below we by and large follow 
the structure that we used in \cite{BKM} for a single Markov modulated infinite-server queue (with a crucial step stemming from \cite{JMTW}), 
and therefore we restrict ourselves to highlighting the main steps. 

\vb

{\it $\circ$ Deviation matrix.} We now introduce some notions that we need across this subsection. For ease, we define them here for the unscaled system, but for the scaled system they can be adapted in a straightforward manner. 

Define by $p_{k\ell}(t):={\mathbb P}({\bs X}(t)=\ell\,|\,{\bs X}(0)=k) = ({\rm e}^{{\bs Q}t})_{k\ell}$ the transition probabilities of ${\bs X}(\cdot).$
The equilibrium probability that block $m$ is `up' is given by $\pi^{(m)} = q_0^{(m)}/(q_0^{(m)}+q_1^{(m)})$, for $m=1,\ldots,K$. 
An important role in the analysis is played by the deviation matrix $D$ (of dimension $\bar K\times \bar K$), whose $(k,\ell)$-th entry is given by
\[D_{k\ell}=\int_0^\infty (p_{k\ell}(t) -\pi_\ell){\rm d}t = \int_0^\infty (({\rm e}^{{\bs Q}t})_{k\ell}-\pi_\ell){\rm d}t,\]
with ${\bs \pi}$, as before, the solution to  ${\bs\pi}{\bs Q} = {\bs 0}$ with entries summing to 1. 

Let $U_k$ be the set of blocks that is `up' when ${\bs X}(t)=k$, and $D_k:=\{1,\ldots,K\}\setminus U_k.$
Define, with $q^{(m)}:= q_0^{(m)}+q_1^{(m)}$,  \[p_{00}^{(m)}(t) :=1-\pi^{(m)}  +\pi^{(m)}  {\rm e}^{-q^{(m)} t},\:\:\:
p_{11}^{(m)}(t) := \pi^{(m)} +(1-\pi^{(m)} ){\rm e}^{-q^{(m)} t}
;\]  in addition $p_{01}^{(m)}(t):=1- p_{00}^{(m)}(t)$ and $p_{10}^{(m)}(t):=1-p_{11}^{(m)}(t)$. Then, as is readily verified,
\begin{equation}
\label{peetje}p_{k\ell}(t)=\left(\prod_{m\in U_k\cap U_\ell} p_{11}^{(m)}(t) \right)\left(\prod_{m\in U_k\cap D_\ell} p_{10}^{(m)}(t)\right)\left(\prod_{m\in D_k\cap U_\ell} p_{01}^{(m)}(t)\right)\left( \prod_{m\in D_k\cap D_\ell} p_{00}^{(m)}(t)\right),\end{equation}
and
\begin{equation}\label{pietje}\pi_\ell =\left( \prod_{m\in U_\ell}\pi^{(m)} \right)\left( \prod_{m\in D_\ell} (1-\pi^{(m)} )\right).\end{equation}
From the explicit expressions for the $p_{ij}^{(m)}(t)$, we conclude that
$p^{(m)}_{j0}(t)-(1-\pi^{(m)})$  and $p^{(m)}_{j1}(t)- \pi^{(m)}$ can be written as a linear combination of terms of the type \[\exp\left(-t\sum_{m\in S} q^{(m)}\right),\] where $S$ is a non-empty set. When subtracting (\ref{pietje}) from (\ref{peetje}), this entails that, for non-empty sets $S_{m'}(k,\ell)$ that depend on $m',$ $k$, and $\ell$,
\[p_{k\ell}(t)-\pi_\ell = \sum_{m'} \alpha_{m'}(k,\ell) \exp\left(-t\sum_{m\in S_{m'}(k,\ell)} q^{(m)}\right),\]
for coefficients $\alpha_{m'}(k,\ell)$ that are straightforward to evaluate but that do not allow an explicit expression. As a consequence, 
\[D_{k\ell} = \sum_{m'} \Bigg(\frac{\alpha_{m'}(k,\ell)}{\displaystyle \sum_{m\in S_{m'}(k,\ell)} q^{(m)}}\Bigg).\]
As an example, we work out the $D$ matrix for the case of one block (i.e.,  $K=1$, or, equivalently, $\bar K=2$). Put $q:=q^{(1)}$, $q_i:=q_i^{(1)}$ (for $i=0,1$), and $\pi:=\pi_1$. Then 
\[D=\int_0^\infty \left(\begin{array}{cc}\pi {\rm e}^{-q t}&-\pi {\rm e}^{-q t}
\\ -(1-\pi){\rm e}^{-q t}
&(1-\pi ){\rm e}^{-q t}\end{array}\right){\rm d}t=
\frac{1}{q^2}
\left(\begin{array}{cc}q_0&-q_0
\\ -q_1
&q_1\end{array}\right)
.
\]
\iffalse
A second example concerns the case in which the $q_i^{(m)}$ do not depend on $m$. With $V_1(k,\ell)$, 
$V_2(k,\ell)$, $V_1(k,\ell)$, and $V_1(k,\ell)$ the number of elements in $U_k\cap U_\ell$, $U_k\cap D_\ell$, $D_k\cap U_\ell$, and $D_k\cap D_\ell$, respectively,
\[p_{k\ell}(t) = \left(p_{11}^{(m)}(t)\right)^{V_1(k,\ell)}\left(p_{10}^{(m)}(t)\right)^{V_2(k,\ell)}\left(p_{01}^{(m)}(t)\right)^{V_3(k,\ell)}\left(p_{00}^{(m)}(t)\right)^{V_4(k,\ell)}.\]\fi

\vb

{\it $\circ$ Time-changed Poisson process representation.}
In this section we repeatedly use the following representation. We throughout use the definition
$\mu_{ijk}:= \mu_{ij} 
{\mathbb I}(i,j,k)$ if $j=1,\ldots,n$ and
\[\mu_{i0k} := \mu_{i0} +  \sum_{j=1,j\not=i}^{n} \mu_{ij} (1- {\mathbb I}(i,j,k)).\]
With all $P_{ij}(\cdot)$ (for $i,j=0,\ldots, n$) independent unit rate Poisson processes, it is directly verified that with the above definition of the rates $\mu_{ijk}$ the numbers of customers in the respective queues  satisfy
\begin{align}\nonumber M^{(N)}_i(t) = P_{0i}(N\lambda_i t) \,+&\, \sum_{j=1,j\not=i}^n P_{ji}\left(\int_0^t M^{(N)}_j(s)\sum_{k=1}^{\bar K}
\mu_{jik} Z_k^{(N)}(s){\rm d}s\right)
\,-\\
&\,\sum_{j=0,j\not=i}^n P_{ij}\left(\int_0^t M^{(N)}_i(s)\sum_{k=1}^{\bar K}
\mu_{ijk} Z_k^{(N)}(s){\rm d}s\right),\label{POIS}\end{align}
with
$Z_k^{(N)}(t) := 1_{\{{\bs X}^{(N)}(t)=k\}}$. This type of Poisson processes with random time-change, and their applications in obtaining scaling-limits,
have been
described in great detail in e.g.\ \cite{AK}.

\vb

{\it $\circ$ SDE for centered and normalized system.}
The idea is to first set up an SDE for ${\bs M}^{(N)}(t)$, which is then translated into
an SDE for its centered and normalized version $\tilde {\bs M}^{(N)}(t)$.

For $i=1,\ldots,n$, by (\ref{POIS}),
\begin{align*}
{\rm d}M^{(N)}_i(t) =&\, N\lambda_i {\rm d}t +
\sum_{j=1,j\not=i}^n \sum_{k=1}^{\bar K} M^{(N)}_j(t) \mu_{jik} Z_k^{(N)}(t)\,{\rm d}t \,-\\
&M_i^{(N)}(t) \sum_{j=0,j\not=i}^n  \sum_{k=1}^{\bar K} 
\mu_{ijk} Z_k^{(N)}(t) \,{\rm d}t+{\rm d}\kappa_i^{(N)}(t),\end{align*}
for some $n$-dimensional martingale ${\bs \kappa}^{(N)}(\cdot).$
In the functional central limit theorem, fluctuations around an average are considered; this $n$-dimensional average vector 
${\bs \varrho}(t)$ solves the following system of linear differential equations:
\[\varrho_i'(t)=\lambda_i + \sum_{j=1,j\not=i}^n \varrho_j(t)\sum_{k=1}^{\bar K}  \mu_{jik} 
\pi_k
-\varrho_i(t) \sum_{j=0,j\not=i}^n  \sum_{k=1}^{\bar K} \mu_{ijk} 
\pi_k ,\]
%\textcolor{red}{SHALL WE ALSO GIVE THE SOLUTION IN MATRIX FORM?}
for $i=1,\ldots,n$. As shown in \cite[Section 3]{JMTW}, this ${\bs \varrho}(\cdot)$ can be considered the fluid limit \cite{WHITT} corresponding to ${\bs M}^{(N)}(\cdot)$. 

Keeping in mind we aim at deriving results for the central-limit regime, we consider
a centered and normalized process whose 
$i$-th component is defined by 
\[\tilde M_i^{(N)}(t) :={N}^{-1/2}\cdot\big(M_i^{(N)}(t) - N\varrho_i(t)\big).\] It is direct that
\begin{align*}
{\rm d}\tilde M^{(N)}_i(t) =&\,   {N}^{1/2}\lambda_i {\rm d}t +
{N}^{-1/2}\sum_{j=1,j\not=i}^n M^{(N)}_j(t)\sum_{k=1}^{\bar K}  \mu_{jik} Z^{(N)}_k(t)\, {\rm d}t \,-\\
&{N}^{-1/2}M_i^{(N)}(t) \sum_{j=0,j\not=i}^n  \sum_{k=1}^{\bar K} 
\mu_{ijk} Z_k^{(N)} (t)\, {\rm d}t-N^{1/2}\varrho'_i(t){\rm d}t+{N}^{-1/2} {\rm d}\kappa_i^{(N)}(t).\end{align*}
The idea now is to plug in the differential equation that is obeyed by ${\bs \varrho}(t)$; the resulting stochastic differential equation resembles the one featuring in \cite{BKM}: omitting a few elementary steps, using the compact notation $\bar Z_k^{(N)}(t):= Z^{(N)}_k(t)-\pi_k$,
\begin{align*}
{\rm d}\tilde M^{(N)}_i(t) =&\,  
N^{-1/2}\sum_{j=1,j\not=i}^n \sum_{k=1}^{\bar K} \mu_{jik}\big(M_j^{(N)}(t)Z^{(N)}_k(t)- N\varrho_j(t) \pi_k\big) {\rm d}t\,-\\
&N^{-1/2}\sum_{j=0,j\not=i}^n \sum_{k=1}^{\bar K} \mu_{ijk}\big(M_i^{(N)}(t)Z^{(N)}_k(t)- N\varrho_i(t) \pi_k\big) {\rm d}t+{N}^{-1/2}
{\rm d}\kappa_i^{(N)}(t)
 \\
=&\,
\sum_{j=1,j\not=i}^n  \sum_{k=1}^{\bar K} \tilde M^{(N)}_j(t)\mu_{jik} Z^{(N)}_k(t)\,{\rm d}t - \sum_{j=0,j\not=i}^n  \sum_{k=1}^{\bar K} \tilde M_i^{(N)}(t)
\mu_{ijk} Z^{(N)}_k(t)\, {\rm d}t\,+\\
&\sqrt{N}\left( \sum_{j=1,j\not=i}^n  \sum_{k=1}^{\bar K}\varrho_j(t) \mu_{jik} \bar Z_k^{(N)}(t)\,{\rm d}t-\sum_{j=0,j\not=i}^n  \sum_{k=1}^{\bar K}\varrho_i(t) \mu_{ijk} \bar Z_k^{(N)}(t)\,{\rm d}t\right)+\\
&{N}^{-1/2}
{\rm d}\kappa_i^{(N)}(t),\end{align*}
or in integral form
\begin{align*}
\tilde M^{(N)}_i(t) =&\,  \int_0^t\sum_{j=1,j\not=i}^n  \sum_{k=1}^{\bar K} \tilde M^{(N)}_j(s)\mu_{jik} Z^{(N)}_k(s)\,{\rm d}s -\int_0^t \sum_{j=0,j\not=i}^n  \sum_{k=1}^{\bar K} \tilde M_i^{(N)}(s)
\mu_{ijk} Z^{(N)}_k(s)\, {\rm d}s\,+\\
&\sqrt{N}\int_0^t\left( \sum_{j=1,j\not=i}^n  \sum_{k=1}^{\bar K}\varrho_j(s) \mu_{jik} \bar Z_k^{(N)}(s){\rm d}s-\sum_{j=0,j\not=i}^n  \sum_{k=1}^{\bar K}\varrho_i(s) \mu_{ijk} \bar Z_k^{(N)}(s){\rm d}s\right)\hspace{-0.9mm}{\rm d}s\,+\\
&{N}^{-1/2}
\kappa_i^{(N)}(t).\end{align*}

$\circ$ {\it A simplification.}
The next step is to verify that as $N\to\infty$, $\tilde M^{(N)}_i(t) -\check M^{(N)}_i(t)$ converges to the zero process as $N\to\infty$ (cf.\ \cite[Lemma 4.4]{JMTW}); here $\check M^{(N)}_i(t)$ is defined as $\tilde M^{(N)}_i(t)$ in the previous display, but now with the $Z_k^{(N)}(t)$ in the first two terms in the right-hand side replaced by $\pi_k.$ To this end, observe that \cite[Thm.\ 5.2]{JMTW}  entails that
\[\int_0^t \tilde M_i^{(N)}(t) \big(Z^{(N)}_k(t)-\pi_k
\big){\rm d}t = \int_0^t \left(\frac{M_i^{(N)}(t)-N\varrho_i(t)}{N}\right) \sqrt{N}\big(Z^{(N)}_k(t)-\pi_k
\big){\rm d}t \to 0,\]
%\textcolor{red}{MAYBE GIVE A REFERENCE TO THEOREM VI 6.22 OF JACOD AND SHIRYAEV?}
using the law-of-large-numbers property that $M_i^{(N)}(t)/N$ converges in probability to $\varrho_i(t)$.  

We have thus arrived at, with $\bar\mu_{ij}:=\sum_{k=1}^{\bar K}\mu_{ijk}\pi_k$, the following system of SDE's:
\begin{align*}
{\rm d}\check M^{(N)}_i(t) =&\, \sum_{j=1,j\not=i}^n   \check M^{(N)}_j(t)\bar\mu_{ji}\,{\rm d}t - \sum_{j=0,j\not=i}^n  \check M_i^{(N)}(t)
\bar \mu_{ij} \, {\rm d}t\,+\\
&\sqrt{N}\left( \sum_{j=1,j\not=i}^n  \sum_{k=1}^{\bar K}\varrho_j(t) \mu_{jik} \bar Z_k^{(N)}(t)\,{\rm d}t-\sum_{j=0,j\not=i}^n  \sum_{k=1}^{\bar K}\varrho_i(t) \mu_{ijk} \bar Z_k^{(N)}(t)\,{\rm d}t\right)+\\
&{N}^{-1/2}
{\rm d}\kappa_i^{(N)}(t).\end{align*}

{\it $\circ$ Functional central limit theorem.}
The following steps echo those in \cite[Section 4]{BKM}. They rely on the transformation
\[{\bs Y}^{(N)}(t) = \exp\big(-{\mathscr M}t  \big) \check {\bs M}^{(N)}(t),\]
with for $i\not=j$ the $(i,j)$-th entry of  the $(n\times n)$-dimensional matrix ${\mathscr M}$ being given by $\bar \mu_{ji}$, whereas the $(i,i)$-th entry is $-\sum_{j=0,j\not=i}^n\bar\mu_{ij}$.
It thus follows that
\[{\rm d}{\bs Y}^{(N)}(t) =  \exp\big(-{\mathscr M}t  \big) \left(\sqrt{N} {\mathscr M}^\circ(t) \bar{\bs Z}^{(N)}(t){\rm d}t+N^{-1/2} {\rm d}{\bs\kappa}^{(N)}(t)\right) ,\]
where the $(i,k)$-th entry of the $(n\times\bar K)$-dimensional matrix ${\mathscr M}^\circ(t)$ is given by 
\[\left({\mathscr M}^\circ(t)\right)_{ik}:=\sum_{j=1,j\not=i}^n  \varrho_j(t) \mu_{jik} -\sum_{j=0,j\not=i}^n \varrho_i(t) \mu_{ijk}.\]
In the next step we analyze the terms $\sqrt{N} {\mathscr M}^\circ(t) \bar{\bs Z}^{(N)}(t)$ and $N^{-1/2} {\rm d}{\bs\kappa}^{(N)}(t)$ separately. 
As in \cite{BKM}, with ${\bs G}^{(N)}(t):=\sqrt{N} {\mathscr M}^\circ(t)\, \bar{\bs Z}^{(N)}(t)$, it follows that ${\bs G}^{(N)}(\cdot) \to 
{\bs G}(\cdot)$ as $N\to\infty$, where ${\bs G}(\cdot)$ satisfies
\[\langle {\bs G}\rangle_t ={\bs V}(t):=\int_0^t {\mathscr M}^\circ(s)\,\Sigma\,({\mathscr M}^\circ(s))^{\rm T}{\rm d}s,\]
with $\Sigma:={\rm diag}\{{\bs \pi}\}D +D^{\rm T}{\rm diag}\{{\bs\pi}\}.$
It entails that ${\bs G}^{(N)}(\cdot)\to {\mathscr M}^\circ(\cdot) \bar {\bs B}(\cdot)$, with $\bar {\bs B}(\cdot)$ a $\bar K$-dimensional zero-mean Brownian motion with covariance matrix $\Sigma$.

Using the precise same argumentation as  in \cite{BKM}, for independent standard Brownian motions $B_{ij}(\cdot)$, with $i,j=0,\ldots,n$,
\[\frac{1}{\sqrt{N} }{\kappa}^{(N)}_i(\cdot) \to 
\sqrt{\lambda_j} B_{0i} (\cdot)+\sum_{j=1,j\not=i}^n\sqrt{\varrho_j(\cdot) \bar\mu_{ji}} B_{ji}(\cdot)-
 \sum_{j=0,j\not=i}^n\sqrt{\varrho_i(\cdot)\bar\mu_{ij}} B_{ij}(\cdot)
,\]
cf.\ (\ref{POIS}); the processes $B_{ij}(\cdot)$ are independent of $\bar {\bs B}(\cdot)$.

Now recall the relation between ${\bs Y}^{(N)}(t)$ and $\check{\bs M} ^{(N)}(t)$, and the fact that    $\tilde {\bs M}^{(N)}(\cdot)-\check {\bs M}^{(N)}(\cdot)$ converges to the zero process as $N\to\infty$. Based on the weak convergence results established above, we thus obtain the following functional central limit theorem. It states that the process under study converges to a (multivariate) process of Ornstein-Uhlenbeck type. 
\begin{proposition}
As $N\to\infty$,  $\tilde {\bs M}^{(N)}(\cdot)$ weakly converges to  $\tilde {\bs M}(\cdot)$,  satisfying the following system of coupled stochastic differential equations:  for $i=1,\ldots,n$,
\begin{align*}
{\rm d}\tilde M_i(t) =&\,\sum_{j=1,j\not=i}^n   \tilde M_j(t)\bar\mu_{ji}\,{\rm d}t - \sum_{j=0,j\not=i}^n  \tilde M_i(t)
\bar \mu_{ij} \, {\rm d}t\,+\\
&\sqrt{\lambda_i}\, {\rm d}B_{0i} (t)+\sum_{j=1,j\not=i}^n\sqrt{\varrho_j(t) \bar\mu_{ji}} {\rm d}B_{ji}(t)-
 \sum_{j=0,j\not=i}^n\sqrt{\varrho_i(t)\bar\mu_{ij}} {\rm d}B_{ij}(t)+({\mathscr M}^\circ(t){\rm d} \bar {\bs B}(t))_i.
\end{align*}
%\textcolor{red}{IN MATRIX NOTATION IT IS CLEARER THAT IT IS AN OU PROCESS}
\end{proposition}
\begin{remark}\label{remlim}
{\em
The distribution of $\tilde {\bs M}(t)$ (for a given value of $t\geqslant 0$, that is) can be explicitly found from known results for multivariate Ornstein-Uhlenbeck processes. If $\tilde {\bs M}(0)$ is constant, then it is an $n$-dimensional Normal distribution with mean ${\bs 0}$ and covariance matrix 
\[ {\mathbb C}{\rm ov}(\tilde{\bs M}(t),\tilde{\bs M}(t))=\int_0^t {\rm e}^{{\mathscr M}(t-s)}{\mathscr M}^\circ(s)\,\Sigma\,({\mathscr M}^\circ(s))^{\rm T}\big({\rm e}^{{\mathscr M}(t-s)}\big)^{\rm T}{\rm d}s + \int_0^t \bar\Sigma(s){\rm d}s,\]
where $\bar\Sigma_{ij}(s)=0$ for $i\not=j$ and
%[not OK yet! -- should be symmetric!]
\[
\bar\Sigma_{ii}(s)=\lambda_i +\sum_{j=1,j\not=i}^n\varrho_j(s)\bar\mu_{ji}-
\sum_{j=0,j\not=i}^n\varrho_i(s)\bar\mu_{ij}.\]
Observe that, as $\bar\Sigma_{ii}(s)=\varrho_i'(s)$ by definition, $\int_0^t \bar\Sigma_{ii}(s){\rm d}s =\varrho_i(t).$ With standard theory on multivariate Ornstein-Uhlenbeck processes also the ${\mathbb C}{\rm ov}(\tilde{\bs M}(t),\tilde{\bs M}(t+u))$, i.e., the  covariance matrix pertaining to the system's time-dependent behavior,  can be determined.
\iffalse
The matrix ${\mathbb C}{\rm ov}(\tilde{\bs M}(t),\tilde{\bs M}(t+u))$, i.e., the  covariance matrix pertaining to the system's time-dependent behavior,  can be characterized along the same lines:
\[{\mathbb C}{\rm ov}(\tilde{\bs M}(t),\tilde{\bs M}(t+u)) = {\mathbb C}{\rm ov}(\tilde{\bs M}(t),\tilde{\bs M}(t))
\,{\rm e}^{{\mathscr M}u},\]
for $t,u\geqslant 0.$ \fi $\hfill\Diamond$}
\end{remark}

\begin{remark}\label{ralpha}{\em As we mentioned, the above result corresponds to the case $\alpha = 1$. Precisely following the line of reasoning of \cite{BKM, JMTW}, for arbitrary $\alpha>0$, we have to define $\tilde M_i^{(N)}(t)$ through
\[\tilde M_i^{(N)}(t) :={N}^{-\beta}\cdot\big(M_i^{(N)}(t) - N\varrho_i(t)\big),\]
with $\beta:= \max\{\tfrac{1}{2}, 1-\tfrac{\alpha}{2}\}.$
 Then the recipe is to go through precisely the same steps as in the proof for $\alpha=1$, but it will turn out that for $\alpha>1$ a specific part of the resulting SDE cancels, whereas for $\alpha<1$ another 
part cancels.

More specifically, it can be argued that for $\alpha>1$ the limiting system of differential equations reduces, for $i=1,\ldots,n$, to
\begin{align*}
{\rm d}\tilde M_i(t) =&\,\sum_{j=1,j\not=i}^n   \tilde M_j(t)\bar\mu_{ji}\,{\rm d}t - \sum_{j=0,j\not=i}^n  \tilde M_i(t)
\bar \mu_{ij} \, {\rm d}t\,+\\
&\sqrt{\lambda_i} {\rm d}B_{0i} (t)+\sum_{j=1,j\not=i}^n\sqrt{\varrho_j(t) \bar\mu_{ji}} {\rm d}B_{ji}(t)-
 \sum_{j=0,j\not=i}^n\sqrt{\varrho_i(t)\bar\mu_{ij}} {\rm d}B_{ij}(t).
\end{align*}
Observe that this entails that  for $\alpha>1$ the limiting system depends on the service rates $\mu_{ijk}$ only through their time averaged counterparts $\bar\mu_{ij}$; this reflects the relatively fast alternating link state process. The system essentially behaves as a network of {\it non-modulated} 
infinite-server queues; in particular, Remark \ref{remlim} indicates that the centered and normalized versions of the individual queues, i.e., the processes $\tilde M_i^{(N)}(\cdot)$, become independent as $N\to\infty$.

For $\alpha\in(0,1)$ the limiting system of differential equations becomes, for $i=1,\ldots,n$,
\begin{align*}
{\rm d}\tilde M_i(t) =&\,\sum_{j=1,j\not=i}^n   \tilde M_j(t)\bar\mu_{ji}\,{\rm d}t - \sum_{j=0,j\not=i}^n  \tilde M_i(t)
\bar \mu_{ij} \, {\rm d}t\,+({\mathscr M}^\circ(t){\rm d} \bar {\bs B}(t))_i.
\end{align*}
In this case, the link state process is relatively slow, such that the scaling limit contains detailed information on the transition rates. In this regime, the individual queues are not asymptotically independent. $\hfill\Diamond$
}\end{remark}

\begin{remark}{\em 
At the expense of some additional notation and administration, the loss process $L^{(N)}(t)$ can be added to vector ${\bs M}^{(N)}(t)$, in that a functional central limit theorem for the centered and normalized version of $(L^{(N)}(t),{\bs M}^{(N)}(t))$ can be established using the same techniques. $\hfill\Diamond$
}\end{remark}

\section{Extensions, ramifications} 
In this section we discuss two extensions. In the first subsection we describe how to adapt the model to incorporate phase-type distributed up- and down-times and phase-type service times. In the second subsection we point out how to adapt the model so as to cover the situation in which blocked customers (i.e., customers wishing to jump from $i$ to $j$ when the link between $i$ and $j$ is down) potentially retry.

\subsection{Phase-type distributions} In case the up- and down-times are of phase-type, this is easily incorporated in that transition rate matrix ${\bs Q}.$ The background process now keeps tracks of each of the links being up or down, but in addition, it gives the phase of the current up- or down-time. If the up-time (down-time, respectively) of the link between $i$ and $j$ is phase type of degree $\delta^{({\rm u})}_{ij}$ ($\delta^{({\rm d})}_{ij},$ respectively),  then the dimension of $X(\cdot)$ is
\[\prod_{i\not =j} (\delta_{ij}^{({\rm u})}+\delta_{ij}^{({\rm d})}).\]
Likewise, the service times can be made phase-type, by keeping track of an infinite-server queue of all clients at any specific node being in a specific phase of the phase-type service time. 

\subsection{Model in which blocked customers retry} 
The previous two section considered the case that all $f_{ij}$ are equal to 1. It is not hard to generalize the results to the situation in which $f_{ij}\in[0,1)$ are allowed as well, but it comes at the price of having to introduce a substantial amount of additional notation. For this reason, we restrict ourselves in the subsection of just pointing out how the results have to be adapted to accommodate $f_{ij}\in[0,1)$. 

It is readily verified that now the joint probability generating function ${\bs  \varphi}(w,{\bs z},t)$ satisfies
\begin{align*}\frac{\partial {\bs  \varphi}(w,{\bs z},t)}{\partial t} = &\,
\sum_{i=1}^n{\bs \varphi} (w,{\bs z},t)\,\lambda_i(z_i-1) \,+\sum_{i=1}^n\sum_{j\not=i}^n\frac{\partial {\bs  \varphi}(w,{\bs z},t)}{\partial z_i}
{\mathbb I}_{\bar K}(i,j)\,
\mu_{ij}\left( {z_j}-{z_i}\right)\,+
 \\
&\sum_{i=1}^n\sum_{j\not=i}^n\frac{\partial {\bs  \varphi}(w,{\bs z},t)}{\partial z_i}{\mathbb J}_{\bar K}(i,j)\,
\mu_{ij}\,f_{ij}\,\left( w-{z_i}\right)\,+\\
&\sum_{i=1}^n \frac{\partial {\bs  \varphi}(w,{\bs z},t)}{\partial z_i}\,\mu_{i0} (1-z_i)
+
 {\bs  \varphi}(w,{\bs z},t)\, {\bs Q}
.\end{align*}
\iffalse

As the derivations are very similar to those for \textcolor{red}{Model I}, we restrict ourselves to presenting the main results. Recall that in this model, no customers are lost.

We start by presenting the partial differential equations characterizing the system's time-dependent behavior. 
\begin{proposition} \label{pr} The joint probability generating function ${\bs  \varphi}({\bs z},t)$ satisfies
\begin{align*}\frac{\partial {\bs  \varphi}({\bs z},t)}{\partial t} = &\,
\sum_{i=1}^n{\bs \varphi} ({\bs z},t)\,\lambda_i(z_i-1) \,+\sum_{i=1}^n\sum_{j\not=i}^n\frac{\partial {\bs  \varphi}({\bs z},t)}{\partial z_i}
{\mathbb I}_{\bar K}(i,j)\,
\mu_{ij}\left( {z_j}-{z_i}\right)\,+
 \\
&\sum_{i=1}^n \frac{\partial {\bs  \varphi}({\bs z},t)}{\partial z_i}\,\mu_{i0} (1-z_i)
+
 {\bs  \varphi}(w,{\bs z},t)\, {\bs Q}
.\end{align*}\fi
The probability generating function of the stationary counterpart ${\bs M}$
follows as before, i.e.,  by plugging in $w=1$ and equating the right-hand side to ${\bs 0}.$

The time-dependent moments can be evaluated from the next result; again, the stationary counterpart follows
by equating the right-hand side to ${\bs 0}$. 
\begin{proposition} \label{P22}
For ${\bs r}\in {\mathbb N}_0^n$, $k\in\{1,\ldots,\bar K\}$ and $t\geqslant 0$,
\begin{align}\nonumber  \frac{\partial \psi_k({\bs r},t)}{\partial t}=\,&\sum_{i=1}^nr_i
\psi_k({\bs r}-{\bs e}_i,t)\lambda_i+ \sum_{i=1}^n\sum_{j=1,j\not=i}^n r_j  \psi_k({\bs r}-{\bs e}_j+{\bs e}_i,t)
\mu_{ijk}^+\,-\\
\,&\sum_{i=1}^n\sum_{j=1,j\not=i}^n r_i \psi_k({{\bs r}},t)\,
(\mu_{ijk}^++f_{ij}\mu_{ijk}^-)- \sum_{i=1}^n  r_i \psi_k({{\bs r}},t)\,
\mu_{i0}+ \sum_{\ell=1}^{\bar K} \psi_\ell({\bs r},t) q_{\ell k}.\label{dve}
\end{align}
\end{proposition}

An interesting special is case is $f_{ij}=0$ for all $i,j$, i.e., the case in which there are no clients lost. 
If in addition full symmetry is assumed, cf.\ Remark \ref{rem1}, the mean allows an explicit expression. It is directly seen that 
for each of the links the mean number of clients present at time $t$, denoted as before by $v(t)$,  satisfies
\begin{equation}\label{fl}v'(t) = \lambda + (n-1)\,v(t) \frac{\nu\pi}{n-1} - v(t)\left((n-1)\, \frac{\nu\pi}{n-1}+\mu_0\right)= \lambda-\mu_0v(t).\end{equation}
It thus follows that
\[v(t) =\frac{\lambda}{\mu_0} \left(1-{\rm e}^{-\mu_0t}\right),\]
which converges to $v:=
{\lambda}/{\mu_0}$
as $t\to\infty$. Observe that $v(t)$ is in this case not affected by $\pi$, due to the fact that parts of the in-flux and out-flux cancel, 
as observed in (\ref{fl}). 

Regarding the functional central limit theorem, the result for $f_{ij}=1$ carries over to that for $f_{ij}\in[0,1)$, but with \[\mu_{i0k} := \mu_{i0} +  \sum_{j=1,j\not=i}^{n} \mu_{ij} f_{ij}(1- {\mathbb I}(i,j,k)).\]  Recall that in this model with a probability $1-f_{ij}$ a customer wishing to jump from node $i$ to node $j$ retries when the link is not available (and hence stays at node $i$).

\section{Examples}
In this section we work out a couple of relevant examples, starting with a tandem network in which the link between the nodes is subject to failure. In the first subsection we consider the case that all blocked customers are lost (i.e., $f_{12}=1$), whereas in the second subsection all blocked customers retry (i.e., $f_{ij}=0$). Section \ref{FCLT_TAND} presents the 
FCLT for the two-node tandem. 
In Section \ref{FCLT_NETW} we consider the FCLT for the case of a symmetric fully connected $n$-node network in which blocked customers are lost (where it is noted that the case in which they retry is dealt with fully analogously); Section \ref{FCLT_RING} deals with its ring-shaped counterpart.

\subsection{Two-node tandem, with blocked customers being lost} We consider a two-node tandem, where traffic arrives at the first node, is sent to the second node after having been served at the first node, and leaves the network after having been served at the second node. Jobs arrive at the first node according to a Poisson process with rate $\lambda$ and have exponentially distributed service times with mean $\mu_i$ at node $i$ ($i=1,2$). The link between node 1 and 2 is up (down, respectively) during an exponentially distributed time with  mean $q_1^{-1}$ ($q_0^{-1}$, respectively). In this subsection we consider the case that $f_{12}=1$: clients who wish to jump from node 1 to node 2 when the link is down, are lost. In this case,
\[{\bs Q} = \left(\begin{array}{rr}-q_0&q_0\\
q_1&-q_1\end{array}\right).\]

Define by $v_{ij}(t)$ the mean number of clients at node $i$ ($i=1,2$) when the background process is in state $j$ ($j=0,1$). Using our expression for the transient first moment, with $p_0(t)$ ($p_1(t)$, respectively) the probability the link is down (up, respectively),
\begin{align*}
v_{10}'(t) &= \lambda p_0(t) + q_1v_{11}(t)  - q_0v_{10}(t)  -\mu_1 v_{10}(t),\\
v_{11}'(t) &= \lambda p_1(t) + q_0v_{10}(t)  - q_1v_{11}(t)  - \mu_1v_{11}(t),\\
v_{20}'(t) &=  q_1v_{21}(t)  - q_0v_{20}(t)  - \mu_2v_{20}(t),\\
v_{21}'(t) &= \mu_1v_{11}(t)+q_0  v_{20}(t) - q_1v_{21}(t)  - \mu_2v_{21}(t).
\end{align*}
This system can be solved in closed form, realizing that the first two differential equations can be solved in isolation first (leading to explicit expressions for $v_{10}(t)$ and $v_{11}(t)$), and then the last two differential equations (using the found expression for $v_{11}(t)$). As these are standard computations involving systems of non-homogeneous linear differential equations, we do not include the expressions here. 

The stationary expectations can be found along the same lines. Let $\Delta_a$ be a (two-dimensional) diagonal matrix, whose diagonal elements are all equal to $a\in{\mathbb R}.$ Then the steady-state means are, with $\pi_0=1-\pi_1 = q_1/(q_0+q_1)$ and ${\bs \pi}=(\pi_0,\:\pi_1)$,
\[(v_{10}, v_{11}) = \lambda\, {\bs\pi} (\Delta_{\mu_1}-{\bs Q})^{-1},\:\:\: (v_{20}, v_{21}) = (0,\:\mu_1v_{11})
(\Delta_{\mu_2}-{\bs Q})^{-1}
.\]
The number of clients lost per unit of time is $\mu_1 v_{10}.$

\subsection{Two-node tandem, with blocked customer retrying} The model in which $f_{12}=0$ has more intricate interactions, as the link between the nodes being down has impact on the number of clients at node 1. It means that in the set of differential equations that we set up for $f_{12}=1$, the first one has to be replaced by $v_{10}'(t) = \lambda p_0(t) + q_1v_{11}(t)  - q_0v_{10}(t).$ Again the time-dependent means can be found in closed form by solving two 2-dimensional systems of linear differential equations, and the stationary means by solving two pairs of linear equations. 
As it turns out, however, we can also explicitly find the distribution of the stationary number of clients residing at node 1, as follows; the resulting formulae reveal the effect of the link failures on the performance experienced at this node. 

Let $\varphi_0(z)$ be the probability generating function of the stationary number of customers at node 1 jointly with the event that the link is down, and $\varphi_1(z)$ its counterpart jointly with the link being up.
The following (differential) equations apply:
\begin{align*}
\lambda(z-1)\varphi_0(z) -q_0 \varphi_0(z)+ q_1\varphi_1(z)&=0,\\
\lambda(z-1)\varphi_1(z) -\mu_1(z-1)\varphi'_1(z) -q_1 \varphi_1(z)+ q_0\varphi_0(z)&=0.
\end{align*}
Inserting $\varphi_0(z)= q_1 \varphi_1(z)/(q_0+\lambda(1-z))$ into the second equation, we obtain
\[ \lambda(z-1)\varphi_1(z) -\mu_1(z-1)\varphi'_1(z) -q_1 \varphi_1(z)+q_0\frac{q_1\varphi_1(z)}{(q_0+\lambda(1-z))}=0,\]
or, equivalently,
\[\frac{\varphi'_1(z)}{\varphi_1(z)} = \frac{\lambda}{\mu_1}\left(1+\frac{q_1}{q_0+\lambda(1-z)}\right).\]
Up to an additive constant, we thus obtain that
\[\log \varphi_1(z) = \frac{\lambda}{\mu_1} z- \frac{q_1}{\mu_1}\log (q_0+\lambda(1-z)),\]
and using that $\varphi_1(1)=\pi_1$, 
\[\varphi_1(z) = \pi_1 \exp\left(\frac{\lambda}{\mu_1}(z-1)\right) \left(\frac{q_0}{q_0+\lambda(1-z)}\right)^{q_1/\mu_1}.\]
Using the relation between $\varphi_0(z)$ and $\varphi_1(z)$, we find that the transform of the stationary number at the first node equals
\[\varphi_0(z)+\varphi_1(z)= \exp\left(\frac{\lambda}{\mu_1}(z-1)\right)
\left(\pi_0\left(\frac{q_0}{q_0+\lambda(1-z)}\right)^{q_1/\mu_1+1}+
\pi_1 \left(\frac{q_0}{q_0+\lambda(1-z)}\right)^{q_1/\mu_1} \right).\]
This expression has the following nice interpretation.
Let $A$ be Poisson with mean $\lambda/\mu_1$, and let $B$ with probability $\pi_0$ be  a negative binomial random variable with parameters $r:=q_1/\mu_1+1$ and $p:= q_0/(q_0+\lambda)$ and with probability $\pi_1$   a negative binomial random variable with parameters $r-1$ and $p$.
Then  the stationary number of customers at the first node is distributed as the sum of two independent random variables $A$ and $B$ (which are both non-negative and integer-valued). Note that if the link would never be down (i.e., $\pi_1=1$ and $q_1=\infty$), the number of customers at node~1 is just Poisson with mean $\lambda/\mu_1$ (like in the ordinary M/M/$\infty$ queue); the additional random variable $B$ (which is a mixture of two negative binomial random variable) thus represents the effect of the link failures. 

\begin{figure}[h]
\includegraphics[height=6.5cm, width=8.1cm]{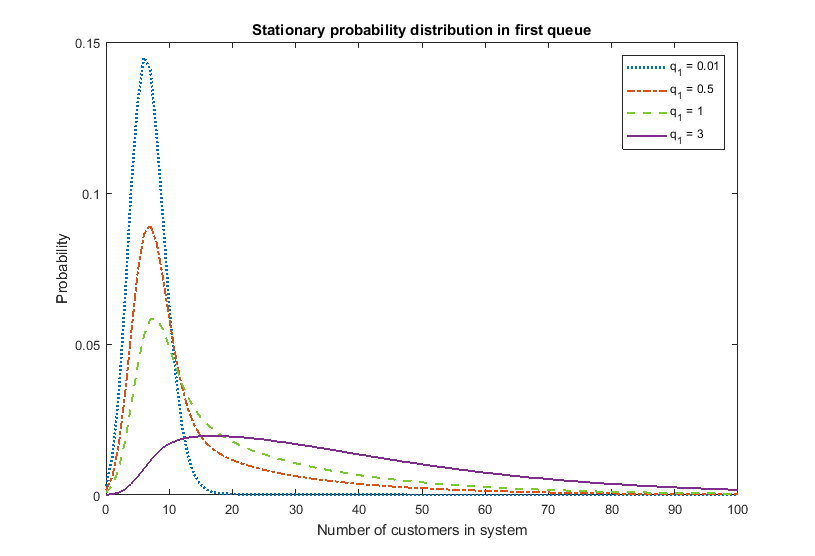}
\includegraphics[height=6.5cm,width=8.1cm]{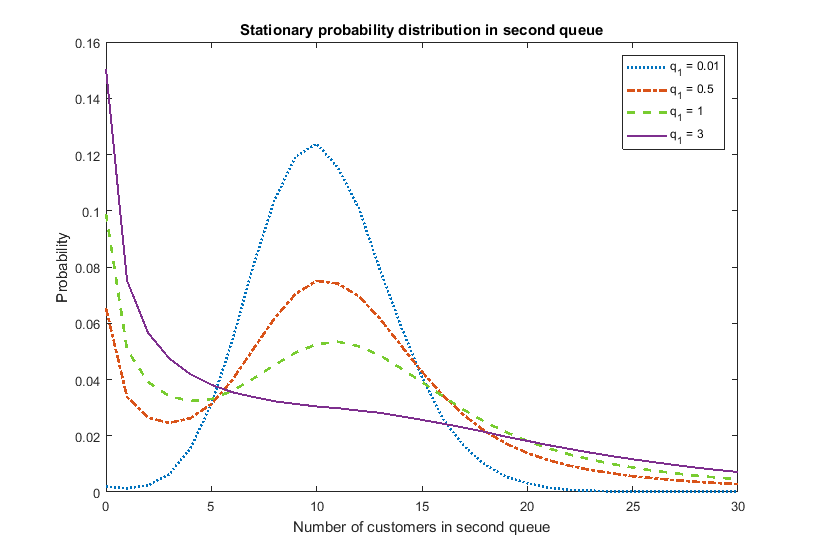}
\caption{\label{P1}Stationary probability density function}
\end{figure}

A numerical illustration is shown in Fig.\ \ref{P1}, where we have estimated 
the stationary random variable ${\bs M}$ by simulation. We have fixed the parameters $\lambda = 20, \mu_1 =3, \mu_2 = 2, q_0=1$ and $f=0$, and  have varied the parameter $q_1$. 
This experiment visualizes  the impact of the link failures on the random variable ${\bs M}$; the left graph corresponds with the upstream queue, and the right graph to the downstream queue. We choose $q_1=0.01$, $0.5$, $1$, and $3$. The simulated numbers in the left panel align with the distribution identified above.  

\subsection{Functional central limit theorem for two-node tandem}  \label{FCLT_TAND}
In this example we derive the FCLT for the two-node tandem. Clients wishing to jump from node 1 to node 2 while the link is down are lost with probability $f:=f_{12}$, and stay at node 1 with probability $1-f.$
First we determine the fluid limit (which is to be used as the `centering function' in our FCLT). We use the same notation as in  the above examples, and in addition we define
$\kappa:=\mu_1(\pi+(1-\pi)f)$ with $\pi:=\pi_1=q_0/(q_0+q_1)$. Then 
\begin{align*}
\varrho'_1(t) &=\lambda -\kappa\varrho_1(t),
\:\:\:
\varrho_2'(t) = \mu_1\pi \varrho_1(t) -\mu_2\varrho_2(t).
\end{align*}
Assuming the system starts empty, we thus obtain
\[\varrho_1(t)=\frac{\lambda}{\kappa}(1-{\rm e}^{-\kappa t}),\:\:\varrho_2(t) =
\frac{\mu_1\pi\lambda}{\kappa}\left(\frac{1-{\rm e}^{-\mu_2 t}}{\mu_2} -\frac{{\rm e}^{-\kappa t}-{\rm e}^{-\mu_2 t}}{\mu_2-\kappa}\right).
\]
It can be checked that
\[{\mathscr M}=  \left(\begin{array}{cc}
-\kappa&0\\
\mu_1\pi&-\mu_2\end{array}\right),\:\:\:
{\mathscr M}^\circ(s) =\left(\begin{array}{cc}
-\varrho_1(s)\mu_1 f&-\varrho_1(s)\mu_1\\
-\varrho_2(s)\mu_2&\varrho_1(s)\mu_1-\varrho_2(s)\mu_2\end{array}\right).\]
In addition, with $q:=q_0+q_1$,
\[\Sigma := \frac{1}{q^2}\left(\begin{array}{cc}
2(1-\pi)q_0 & -(1-\pi)q_0 - \pi q_1\\
-(1-\pi)q_0 - \pi q_1 & 2\pi q_1\end{array}\right).\]
With 
\[{\rm e}^{{\mathscr M}t} = \left(\begin{array}{cc}
{\rm e}^{-\kappa t}&0\\
\mu_1\pi({\rm e}^{-\mu_2 t}-{\rm e}^{-\kappa t})/(\kappa-\mu_2)&{\rm e}^{-\mu_2 t}\end{array}\right),\]
the matrix ${\mathbb C}{\rm ov}(\tilde{\bs M}(t),\tilde{\bs M}(t))$ can be evaluated using the expressions presented in Remark \ref{remlim}.

We now present a number of figures that illustrate the applicability 
of the diffusion limit as an approximation to the original population process. First we consider the scaling parameter $N=100$ and the   parameters  $\lambda=25$, $\mu_1=10$, $\mu_2=20,$ and $f=0$. In the two simulations we varied the transition rates: they are $q_0=0.2$, $q_1=0.6$ in the left sample path, and $q_0=30$, $q_1=20$ in the right sample path.
We pick $\alpha=1$, so that the arrival rate $\lambda^{(N)} $ is $N\lambda$, whereas transition rates are set to $q_0^{(N)} = N q_0$ and  $q_1^{(N)}=Nq_1$. 
\begin{figure}[h]
\includegraphics[height=6.5cm, width=7.8cm]{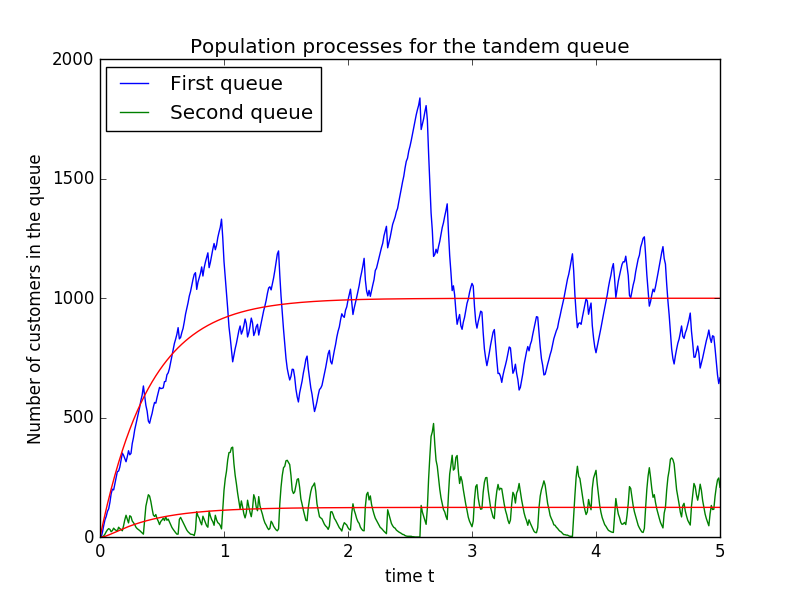}
\includegraphics[height=6.5cm,width=7.8cm]{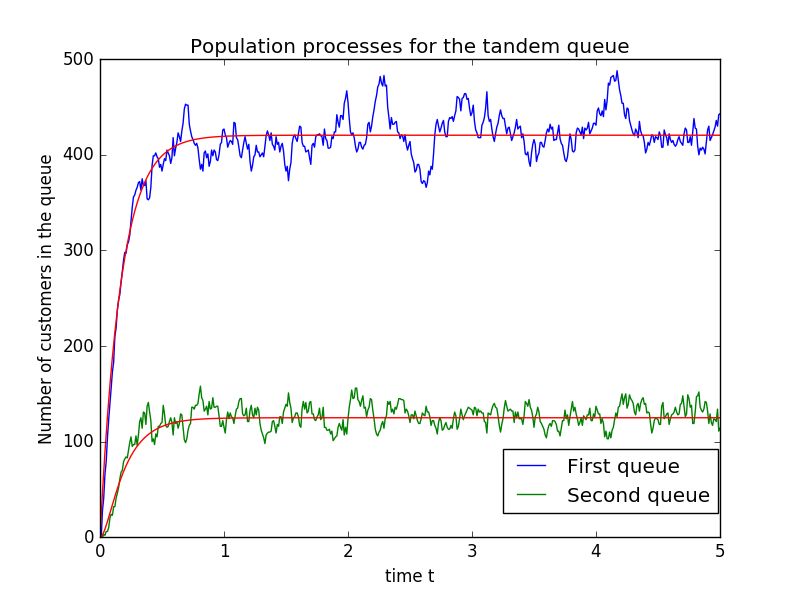}
\caption{\label{P2}Sample paths and centering functions.}
\end{figure}
The red curves appearing in Fig.\ \ref{P2} above correspond to the functions $N\varrho_1(\cdot)$ and $N\varrho_2(\cdot)$, where $\varrho_1(\cdot)$ and $\varrho_2(\cdot)$ are the two 
`centering functions' that were computed above. The blue and green curves are the corresponding sample paths.

In Fig.\ 3  histograms are presented for the centered and scaled 
population process in each queue. The parameters $\lambda$, $\mu_1$, $\mu_2$ and $f$ are chosen as above; the transition rates of the background process are  $q_0=30$ and $q_1=20$. We consider $t=2$ and $N=60$. The graphs show that the Gaussian limiting distribution provides an accurate fit; the dotted curves correspond to the zero-mean Gaussian distribution with the variance in line with the diffusion limit. 

\begin{figure}[h]
\includegraphics[height=6.5cm, width=7.8cm]{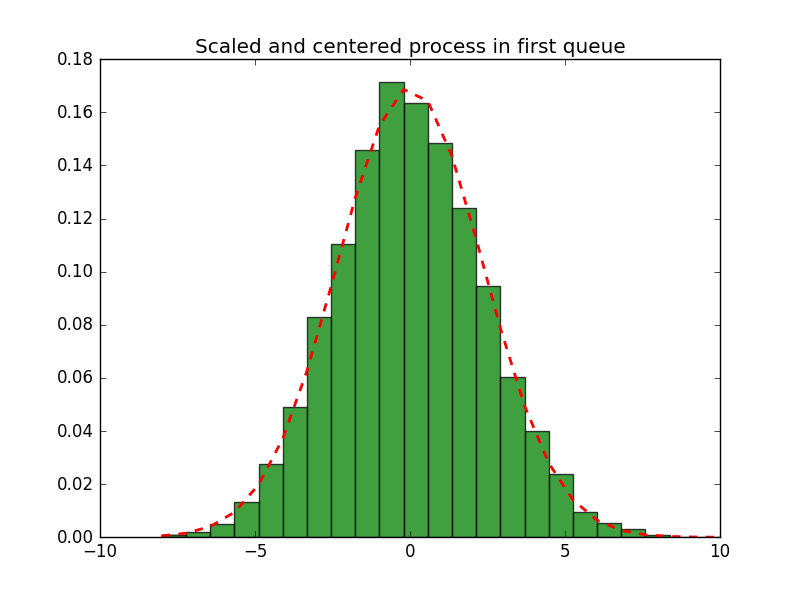}
\includegraphics[height=6.5cm,width=7.8cm]{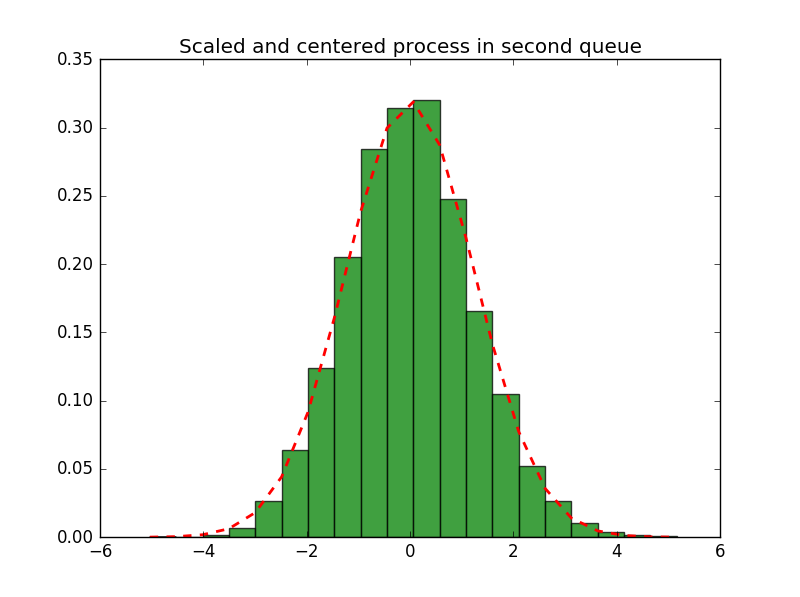}
\caption{\label{P3} Histograms for centered and scaled stationary population.}
\end{figure}

\subsection{Functional central limit theorem for symmetric fully connected one-block network}  \label{FCLT_NETW}In this subsection we consider the functional central limit theorem for a network with just a single block (i.e., all links alternate between being `up' and `down' simultaneously), and all parameters chosen symmetrically; blocked customers are assumed lost (but the case in which they retry works analogously). 

More concretely, the situation considered is the following. The arrival rate at each node is $N\lambda$. The down-time of all links is exponentially distributed with mean $(Nq_0)^{-1}$, whereas the up-time is exponentially distributed with mean $(Nq_1)^{-1}$. As in Remark \ref{rem1}, the service rate is $\sigma:=\nu+\mu_0$; after service completion a  client leaves the network with probability $\mu_0/\sigma$ and wants to move to another node (picked uniformly at random) with probability $\nu/\sigma.$ Define $\pi:=q_0/(q_0+q_1).$

First we find the `centering function' $\varrho(\cdot)$:
\[\varrho'(t) = \lambda + (n-1)\varrho(t) \frac{\nu\pi}{n-1} - \varrho(t)\sigma, \]
solved by, assuming the queues start empty and defining $\kappa:=\nu(1-\pi)+\mu_0$,
\[\varrho(t) =\frac{\lambda}{\kappa}\left(1-{\rm e}^{-\kappa t}\right).\]
Recalling that $\bar\mu_{ij}:=\sum_{k=1}^{\bar K}\mu_{ijk}\pi_k$, for $i\not=j$ and $i,j\in\{1,\ldots,n\}$, we find that ${\mathscr{M}(i,j)=\bar\mu_{ji}} =  \nu\pi/(n-1)$. In addition, $\mathscr{M}(i,i)=\bar{\mu}_{i0}=-\sigma$. As a consequence, with $E_n$ an $n\times n$ all-ones matrix,
\[{\mathscr M}=\omega_1 E_n+ \omega_2I_n,\:\:\:\:
\omega_1:= \frac{\nu\pi}{n-1},\:\:\:\omega_2:=-\left(\frac{\nu\pi}{n-1}+\sigma\right) .\]
It is readily checked that
\[{\rm e}^{{\mathscr M}t} =  \left(E_n\left(\frac{{\rm e}^{\omega_1nt}-1}{n}\right)+I_n\right)\,{\rm e}^{\omega_2t}.\]
The matrix ${\mathscr M}^\circ(s)$ is an $(n\times 2)$-dimensional matrix whose entries in the first column (which are corresponding to the links being down) are all $m_0(s):= -\varrho(s)\sigma$, and whose entries in the second column (which are corresponding to the links being up) are all $m_1(s):= -\varrho(s)\mu_0$. 
The matrix $\Sigma$ is as in Section \ref{FCLT_TAND}.

Using the expressions from Remark \ref{remlim}, we obtain, with ${\bs 1}_n$ the $n$-dimensional all-ones vector and ${\bs m}(s):=(m_0(s),m_1(s))^{\rm T}$,
\begin{equation}
\label{COV} {\mathbb C}{\rm ov}(\tilde{\bs M}(t),\tilde{\bs M}(t))=\int_0^t {\rm e}^{{\mathscr M}(t-s)}{\bs 1}_n ({\bs m}(s))^{\rm T}\,\Sigma\,{\bs m}(s)\,{\bs 1}_n^{\rm T}\,\big({\rm e}^{{\mathscr M}(t-s)}\big)^{\rm T}{\rm d}s + {\rm diag}\{{\bs \varrho}(t)\}.\end{equation}
The next step is to explicitly evaluate the integral. 
To this end, we first observe that
${\rm e}^{{\mathscr M}t}\,{\bs 1}_n ={\rm e}^{(\omega_1 n+\omega_2)t} \,{\bs 1}_n={\rm e}^{-\kappa t}\,{\bs 1}_n.$

We now evaluate $({\bs m}(s))^{\rm T}\,\Sigma\,{\bs m}(s)$. 
With ${\bs x}\in{\mathbb R}^{\bar K}$, $\alpha\in {\mathbb R}$, and ${\bs 1}\equiv {\bs 1}_{\bar K}$, 
\begin{align*}
({\bs x}- \alpha {\bs 1})^{\rm T}\Sigma ({\bs x}- \alpha {\bs 1})&= ({\bs x}- \alpha {\bs 1})^{\rm T}
({\rm diag}\{\bs\pi\}D + D^{\rm T} {\rm diag}\{\bs\pi\}) ({\bs x}- \alpha {\bs 1})\\
&=  ({\bs x}- \alpha {\bs 1})^{\rm T}
{\rm diag}\{\bs\pi\}D {\bs x} + {\bs x}^{\rm T} D^{\rm T} {\rm diag}\{\bs\pi\} ({\bs x}- \alpha {\bs 1})
\end{align*}
due to $D{\bs 1}={\bs 0}.$ In addition, ${\bs \pi}^{\rm T}D={\bs 1}^{\rm T} {\rm diag}\{\bs\pi\} D={\bs 0}$, so that
\[({\bs x}- \alpha {\bs 1})^{\rm T}\Sigma ({\bs x}- \alpha {\bs 1})=
  {\bs x}^{\rm T}
{\rm diag}\{\bs\pi\}D {\bs x} + {\bs x}^{\rm T} D^{\rm T} {\rm diag}\{\bs\pi\}  {\bs x}=
 {\bs x}^{\rm T}\Sigma  {\bs x}.\]
We therefore have that, when evaluating $({\bs m}(s))^{\rm T}\,\Sigma\,{\bs m}(s)$, we can replace ${\bs m}(s)$ by ${\bs m}(s) +\varrho(s)\mu_0\,{\bs 1}$, and consequently
\begin{align*}({\bs m}(s))^{\rm T}\,\Sigma\,{\bs m}(s) &= q^{-2}{(-\varrho(s)\nu,\,0)}\left(\begin{array}{cc}
2(1-\pi)q_0 & -(1-\pi)q_0 - \pi q_1\\
-(1-\pi)q_0 - \pi q_1 & 2\pi q_1\end{array}\right){\left(\begin{array}{c}
-\varrho(s)\nu\\0\end{array}\right)}\\
&= 2 q^{-2} \,(\varrho(s))^2\,\nu^2(1-\pi)q_0 =2q_0q_1\,(\varrho(s))^2\,\nu^2/q^3.\end{align*}

Noting that ${\bs 1}_n\,{\bs 1}_n^{\rm T}=E_n$, we conclude that (\ref{COV}) can be written as 
$\xi_n(t) E_n+{\rm diag}\{{\bs \varrho}(t)\}$, where
\begin{align*}\xi_n(t):=&\,\,2 q_0q_1\frac{\lambda^2\nu^2}{\kappa^2q^3} \, \int_0^t (1-2{\rm e}^{-\kappa s}+{\rm e}^{-2\kappa s}) {\rm e}^{-2\kappa(t-s)}{\rm d}s\\
=&\,\,2 q_0q_1 \frac{\lambda^2\nu^2}{\kappa^2q^3}  \left(\frac{1-{\rm e}^{-2\kappa t}}{2\kappa}
-2{\rm e}^{-\kappa t}\frac{1-{\rm e}^{-\kappa t}}{\kappa}+t\,{\rm e}^{-2\kappa t}\right).
\end{align*}
With $\varrho:=\lambda/\kappa$, we also obtain
\[\lim_{t\to\infty} {\mathbb C}{\rm ov}(\tilde{\bs M}(t),\tilde{\bs M}(t))=
q_0q_1 \frac{\lambda^2\nu^2}{\kappa^3q^3} \,E_n+ {\rm diag}\{{\bs \varrho}\}.\]
%Observe that in the limiting regime (i.e., $N\to\infty$) the variances and covariances scale linearly in the number of queues $n$.

\subsection{Functional central limit theorem for symmetric ring-shaped one-block network}  \label{FCLT_RING}
The setting considered is the same as in the previous subsection, with the only exception that a job served at queue $m$ moves to queue $m+1$ (where $n+1$ is to be understood as $1$). More specifically, the service rate is $\sigma:=\nu+\mu_0$; after service completion a  client leaves the network with probability $\mu_0/\sigma$ and wants to move to the next node with probability $\nu/\sigma.$ We concentrate on the case that $f_{i,i+1} = f_{n,1}=1$ (i.e., during outages jobs that wish to jump to the next queue are lost), but we remark that the case of retry can be handled analogously. The shape of the centering function $\varrho(\cdot)$ is as in the previous subsection, with the same $\kappa=\nu(1-\pi)+\mu_0.$

It is verified that $\bar\mu_{i,i+1}=\bar\mu_{n,1}=\nu\pi$ (for $i=1,\ldots n-1$) and $\bar\mu_{ii}=-(\nu\pi+\sigma)$. As a consequence, with $F_n$ denoting an $n\times n$ matrix with ones on the subdiagonal and at entry $(1,n)$,
\[{\mathscr M}=\omega_1 F_n +\omega_2 I_n,\:\:\:\omega_1:=\nu\pi,\:\:\;\omega_2=-(\nu\pi+\sigma).\]
The matrix ${\mathscr M}^\circ(s)$ is, analogously to what we found in Section \ref{FCLT_NETW}, an $(n\times 2)$-dimensional matrix whose entries in the first column are all $m_0(s):= -\varrho(s)\sigma$, and whose entries in the second column are all $m_1(s):= -\varrho(s)\mu_0$. 
The matrix $\Sigma$ is as defined in Section \ref{FCLT_TAND}.

Observe that $F_n{\bs 1}_n={\bs 1}_n$ (as $F_n$ is a permutation matrix), and therefore $F_n^{\,k}{\bs 1}_n={\bs 1}_n$, so that
\[{\rm e}^{\omega_1 F_n t} {\bs 1}_n = \sum_{k=0}^\infty \frac{F_n^{\,k}{\bs 1}_n}{k!}(\omega_1 t)^k =
{\rm e}^{\omega_1 t}{\bs 1}_n.\]
This allows us to conclude that ${\rm e}^{{\mathscr M}t} {\bs 1}_n={\rm e}^{-\sigma t}{\bs 1}_n$. The remaining computations are as in the previous subsection. 
%This entails that all non-diagonal elements of ${\mathbb C}{\rm ov}(\tilde{\bs M}(t),\tilde{\bs M}(t))$ are identical. 

\section{Concluding remarks} In this paper we have considered networks of infinite-server queues with faulty links. Clients that wish to jump from one queue to another while the required link is down are with a given probability lost (and otherwise stay at the origin node to retry after an exponentially distributed amount of time). For this model we derived prelimit results (in terms of differential equations uniquely characterizing the probability generating function, as well as a recursion by which all moments can be determined)  as well as a functional central limit theorem (after appropriately scaling the arrival rates and the links' failure and repair rates). 

This work is among the first papers on queueing processes on dynamically evolving  random graphs. Several alternative models can be considered; we mention a few here. (i) In our work all queues were of infinite-server type. In many applications, one would rather be interested in the queueing discipline being single-server of many-server. (ii) Our probabilistic analysis covers means and diffusion limits, but extreme behavior (`far away from the mean') is not included. Such a rare-event analysis sheds light on the probability that the queueing process attains values in remote sets. (iii) In dynamically evolving networks, typically measures are taken when links fail; think of rerouting mechanisms. This makes the systematic study of the efficacy of such rerouting protocols a relevant topic for further study.

{\small
\section{Acknowledgments}
The authors thank  Melike Baykal-G\"ursoy (Rutgers University, USA), Dieter Fiems (University of Ghent, Belgium), Brendan Patch (The University of Queensland, Australia \& University of Amsterdam, the Netherlands), and Peter Taylor (The University of Melbourne, Australia) for helpful comments and discussions. 
}

\bibliographystyle{plain}

{\small \begin{thebibliography}{111}

%\bibitem{KS}
%{\sc O. Kella} and {\sc W. Stadje} (2002). Markov modulated linear fluid networks with Markov additive input. {\it Journal of Applied Probability} {\bf 39}, pp.\ 413-420.



%\bibitem{KW}
%{\sc O. Kella} and {\sc W. Whitt} (1999). 
%Linear stochastic fluid networks.
%{\it Journal of Applied Probability} {\bf 36}, pp.\ 244-260.
\bibitem{AND}
{\sc D. Anderson, J. Blom, M. Mandjes, H. Thorsdottir,} and {\sc K. De Turck} (2016). A functional central limit theorem for a Markov-modulated infinite-server queue. {\it Methodology and Computing in Applied Probability} {\bf 18}, pp.\  153-168.
\bibitem{AK}
{\sc D. Anderson} and {\sc T. Kurtz} (2011). Continuous-time Markov chain models for chemical reaction networks. In: H. Koeppl, G. Setti, M. di Bernardo, D. Densmore  (eds), {\it Design and Analysis of Biomolecular Circuits}. Springer.
%\bibitem{JA} {\sc J. Artalejo} and {\sc A. G\'omez-Corral} (2008). {\it Retrial Queueing Systems: a Computational Approach}.  Springer. 
%\bibitem{AsmussenAPQ}
%{\sc S.~Asmussen} (2004). {\em Applied Probability and Queues}. Springer.
\bibitem{AV}
{\sc L. Avena, H. Guldas, R.van der Hofstad,} and {\sc F. den Hollander} (2017).
Mixing times of the non-backtracking random walk on dynamic configuration models.
 {\it ArXiv}: {\tt \footnotesize 1606.07639}.
\bibitem{KUR}
{\sc
K. Ball, T. Kurtz, L. Popovic,} and {\sc G. Rempala} (2006). Asymptotic analysis of multiscale approximations to reaction networks. {\it Annals of Applied Probability} {\bf 16}, pp.\ 1925-1961.
\bibitem{BKM}
{\sc J. Blom, K. De Turck,} and {\sc M. Mandjes} (2016). Functional central limit theorems for Markov-modulated infinite-server systems. {\it Mathematical Methods of Operations Research} {\bf 83}, pp. 351-372.
\bibitem{BKMT}
{\sc J. Blom, O. Kella, M. Mandjes}, and {\sc H. Thorsdottir} (2014). Markov-modulated infinite-server queues with general service times. {\it Queueing Systems} {\bf 76}, pp. 403-424.
\bibitem{DAU}
{\sc B. D'Auria}  (2008). M/M/$\infty$ queues in semi-Markovian random environment. {\it Queueing Systems} {\bf 58}, pp. 221-237.
\bibitem{GUR}
{\sc M. Baykal-G\"ursoy} and {\sc W. Xiao} (2004).
Stochastic decomposition in M/M/$\infty$ queues with Markov modulated service rates. {\it Queueing Systems} {\bf 48}, pp.\ 75-88.
\bibitem{BM}
{\sc I. Benjamini} and {\sc E. Mossel} (2003). On the mixing time of a simple random walk on
the supercritical percolation cluster. {\it Probability Theory and Related Fields}
{\bf 125}, pp.\ 408-420.
\bibitem{BR1}
{\sc N. Berestycki, E. Lubetzky, Y. Peres,} and {\sc A. Sly} (2015). Random walks on the random
graph.  {\it ArXiv}: {\tt \footnotesize 1504.01999}.
\bibitem{DH}
{\sc S. Dharmaraja, A. Di Crescenzo, V. Giorno,} and {\sc A. Nobile} (2015). A continuous-time Ehrenfest model with catastrophes and its jump-diffusion approximation. {\it Journal of Statistical Physics} {\bf 161}, pp. 326-345.
\bibitem{ER}
{\sc P. Erd\H{o}s} and {\sc A. R\'enyi} (1959). On random graphs I. {\it Publicationes Mathematicae Debrecen} {\bf 6}, pp. 290-297.
\bibitem{FMP}
{\sc D. Fiems, M. Mandjes}, and {\sc B. Patch} (2017). A network of infinite-server queues with multiplicative transitions. {\it Preprint.}
\bibitem{FRA}
{\sc B. Fralix} and {\sc I.  Adan} (2009). An infinite-server queue influenced by a semi-Markovian environment. {\it Queueing Systems} {\bf 61}, pp.\ 65-84.
\bibitem{EG}
{\sc E. Gilbert} (1959). Random Graphs. {\it Annals of Mathematical Statistics} {\bf 30}, pp.\ 1141-1144.
\bibitem{HO1}
{\sc P. Holme} and  {\sc J. Saram\"aki} (2012). Temporal networks. {\it Physics Reports} {\bf 519}, pp.\ 97-125.
\bibitem{HO2}
{\sc P. Holme} (2015). Modern temporal network theory: a colloquium. {\it European Physical Journal B} {\bf  88}, pp. 1-30.
%\bibitem{MAR} {\sc H.M. Jansen} (2018).
%Scaling limits for modulated infinite-server queues and related stochastic processes. Ph.D. Thesis, University of Ghent and University of Amsterdam. 
\bibitem{JMTW}
{\sc H.M. Jansen, M. Mandjes, K. De Turck,} and {\sc S. Wittevrongel} (2017).
Diffusion limits for networks of {M}arkov-modulated  infinite-server queues. {\it ArXiv}: {\tt \footnotesize 1712.04251}. 
\bibitem{KEL}
{\sc F. Kelly} (1979). 
{\it Reversibility and Stochastic Networks.}  Wiley.
\bibitem{KEI}
{\sc J. Keilson} and {\sc L. Servi} (1993). The matrix M/M/$\infty$
  system: retrial models and Markov modulated sources. {\it Advances in Applied Probability} {\bf 25}, pp.\ 453-471.
\bibitem{KR}
{\sc J. Kurose} and {\sc K. Ross} (2004). {\it Computer Networking}, 3rd ed. Benjamin/Cummings.
\bibitem{LOU}
{\sc G. Louchard} (1988).
Large finite population queueing systems part I: the infinite-server model.
{\it Stochastic Models} {\bf 4}, pp.\ 473-505.
\bibitem{MSBS}
{\sc 
M. Mandjes, N.J. Starreveld, R. Bekker}, and {\sc P. Spreij} (2018).
Dynamic Erd\H{o}s-R\'enyi graphs. {\it Lecture Notes in Computer Science} {\bf 10000}, to appear. {\it ArXiv}: {\tt \footnotesize 	1703.05505}. 
\bibitem{KOEN}
{\sc M.  Mandjes} and {\sc K. De Turck} (2016). Markov-modulated infinite-server queues driven by a common background process. 
{\it Stochastic Models} {\bf 32}, pp.\ 206-232. 
\bibitem{MW}
{\sc W. Massey} and {\sc W. Whitt} (1993).
Networks of infinite-server queues with nonstationary Poisson input.
{\it Queueing Systems} {\bf 13},  pp.\ 183-250.
\bibitem{OCP}
{\sc C. O'Cinneide} and {\sc P. Purdue} (1986) The M/M/$\infty$ queue in a random environment. {\it Journal of Applied Probability} {\bf 23}, pp.\ 175-184.
\bibitem{BRUG}
{\sc A. Schwabe, M. Dobrzy\'ski, K. Rybakova, P. Verschure,} and {\sc F. Bruggeman} (2011). Origins of stochastic intracellular processes and consequences for cell-to-cell variability and cellular survival strategies. {\it Methods in Enzymology}, {\bf 500}, pp.\ 597-625.
\bibitem{SER}
{\sc R. Serfozo} (1999). {\em Introduction to Stochastic Networks}. Springer.
\bibitem{ZMN} {\sc X. Zhang, C. Moore,} and {\sc M. Newman} (2017). Random graph models for dynamic
networks. {\it European Physical Journal B}, {\bf 90},  200.   {\it ArXiv}: {\tt \footnotesize 1607.07570v1}.
\bibitem{WHITT} {\sc W. Whitt} (2002).  {\it Stochastic-Process Limits}. Springer.


\end{thebibliography}}

\end{document}